\documentclass[12pt]{elsarticle}
\usepackage{hyperref}
\usepackage{amsmath}
\usepackage{graphicx}
\usepackage{graphics}
\usepackage{wrapfig}
\usepackage{epsfig}
\usepackage{upgreek}
\usepackage{mathtools}
\usepackage{epstopdf}
\usepackage{amsmath}
\usepackage{bigints}
\usepackage{times}
\usepackage{amssymb}
\usepackage{array}
\usepackage{latexsym}
\usepackage{amsthm}
\usepackage{amsfonts}
\usepackage{nomencl}
\usepackage{lineno}

\makenomenclature
\usepackage{graphics}
\usepackage{hyperref}
\usepackage{array,tabularx}
\usepackage{scalerel,stackengine}
\usepackage[english]{babel}
\allowdisplaybreaks

\newtheorem{remark}{Remark}

\stackMath
\numberwithin{equation}{section}
\date{}
\newtheorem{t1}{Theorem}[section]

\newtheorem{l1}{Lemma}[section]

\usepackage[utf8]{inputenc}
\allowdisplaybreaks

\begin{document}
\begin{frontmatter} 
\title{On a degenerate boundary value problem to relativistic magnetohydrodynamics with a general pressure law}
\author{Rahul Barthwal}
\author{T. Raja Sekhar}
\address{Department of Mathematics, Indian Institute of Technology Kharagpur, Kharagpur,  India} 
{
\begin{abstract} 
{This work is concerned with establishing the existence and uniqueness of the solution to a mixed-type degenerate boundary value problem for a relativistic magnetohydrodynamics system. We first consider a full relativistic magnetohydrodynamics system and reduce it to a simplified form under the assumption that the magnetic field vector is orthogonal to the velocity vector. We consider a boundary value problem for the steady part of this reduced system where the boundary data is prescribed on a sonic boundary and a  characteristic curve. Here the main difficulty is the consideration of a relativistic system, with a general equation of state while considering the magnetic field effects as well, which we believe has never been analyzed before in the context of the analytical study of sonic-supersonic flows. Also, the degeneracy of the governing equations along the sonic curve is a crucial challenge. However, we employ the iteration method used in the work of Li and Hu  \cite{li2019degenerate} to prove the existence and uniqueness of a local classical supersonic solution in the partial hodograph plane first and finally, we recover a local smooth solution to the boundary value problem in the physical plane by applying an inverse transformation.}
\end{abstract}
}
\end{frontmatter} 
{
\section{Introduction}
One of the most significant issues in the mathematical analysis of gas dynamical systems is the analysis of transonic flow which is present in a variety of significant physical phenomena. The analysis of sonic-supersonic patches is crucial in the mathematical study of transonic flow problems. According to Courant and Friedrichs \cite{courant0}, if the Mach number of the flow passing through an obstacle or duct is not much below one, the flow can become supersonic near the duct's surface due to the convexity of the duct and reverts back to being purely subsonic. Many scientific and aerospace applications naturally include similar circumstances, examples including the flow via a nozzle or a flow over a symmetric or asymmetric airfoil. A global solution to these transonic flow problems has been the subject of numerous substantial contributions over the past century, but it is still an open mathematical problem even now. Because a transonic structure is made up of supersonic and subsonic components that are separated by either a transonic shock or a sonic boundary, transonic flows are much more complex to analyze in contrast to a purely supersonic or subsonic flow. Also, the governing equations of the flow are generally linearly degenerate on a sonic boundary and may change their nature from hyperbolic (supersonic) to elliptic (subsonic) across the degenerate boundary; see \cite{li1998two, li2009interaction, li2011characteristic, barthwal2022existence, barthwal2021existence, barthwal2022simple}. 

According to Morawetz's research on transonic flow in channels and ducts (see Figure \ref{fig: 1}) a smooth flow is absent in general. This suggests that in the downstream flow, there could be a transonic shock  \cite{morawetz1964non}. The nonexistence of a smooth transonic flow around a plane blunt airfoil was proved by Morawetz \cite{morawetz1957non} using the compensated-compactness framework. The same framework was utilized to prove the existence of weak solutions to transonic flow problems by Morawetz \cite{morawetz1954uniqueness} and Chen et al. \cite{chen2007two}. Recent years have seen the development of much significant research for the subsonic-sonic and sonic-supersonic parts of the transonic flow for two-dimensional gas dynamics systems and other related mixed-type systems. For the subsonic part, the global existence of solutions has been developed extensively in the last few decades. For a subsonic-sonic part of a long nozzle, the global existence of a solution was proved by Xie and Xin \cite{xie2007global, xie2010global}. Further, they proved the well-posedness of the sonic-subsonic and subsonic flows with critical mass flux for isentropic Euler equations. Their work was extended to the full Euler equations by Chen et al. \cite{chen2012global}. Chen et al. \cite{chen2016subsonic} developed the sonic-subsonic limit of approximate solutions for the multi-dimensional full Euler equations. Recently, Wang and Xin studied Meyer-type transonic flows in de Laval nozzles and established the existence and uniqueness of a smooth solution. Also, they proved the well-posedness of global subsonic-sonic flow problems \cite{wang2021regular}. As far as the sonic-supersonic part is concerned, the relevant results are still very limited. Here we give a brief overview of some significant developments in the context of supersonic-sonic flows over the last few years. Under the assumption that the flow is irrotational, a local sonic-supersonic solution was developed by Zhang and Zheng \cite{zhang2014sonic} for isentropic Euler equations in two dimensions. Their results were extended to the case of two-dimensional steady and unsteady full Euler equations by Hu and Li \cite{hu2019sonic, hu2020sonic}. Further, they also developed a global smooth sonic-supersonic solution and discussed its behaviour near the sonic boundary in \cite{hu2020global}. When tackling sonic-supersonic boundary value problems, the partial hodograph mapping employed in the work of Li and Hu has been proven to be extremely crucial; see \cite{li2019degenerate, hu2019sonic, hu2020global} and references cited therein for more details. We also refer interested readers to \cite{du2011subsonic, du2014steady, chen2016two, wang2013degenerate, hu2021sonic, li2019degenerate, chen2007two} for more such results on subsonic-sonic and sonic-supersonic patch problems.

\begin{figure}
    \centering
    \includegraphics[width= 4 in]{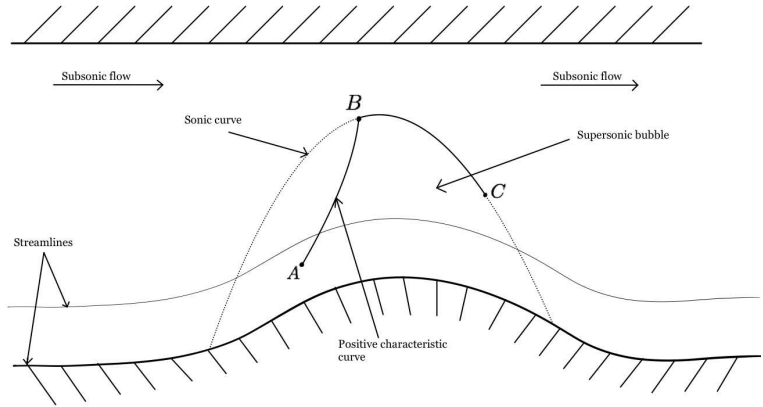}
    \caption{Transonic flow in a duct with a supersonic bubble}
    \label{fig: 1}
\end{figure}

The velocity of fluid particles in the stellar to galactic scales (active galactic nuclei, core-collapse supernovae, superluminal jets, coalescing neutron stars, gamma-ray bursts and the formation of black holes etc.) in the areas of plasma physics, astrophysics, and nuclear physics is typically very high and frequently very close to the speed of light, necessitating the consideration of relativistic effects and rendering obsolete the classical Euler equations. The governing system in question for such a high-speed flow is known as the relativistic Euler system. Recently, a number of sonic-supersonic boundary value problems for a two-dimensional relativistic Euler system and axisymmetric relativistic Euler system has been considered; see viz. \cite{fan2022sonic, barthwal2023existence}. However, in the analytical study of the mathematical theory of transonic flows, magnetic effects for such complicated systems have not been considered before to the best of the authors' knowledge. Throughout celestial objects as well as the entire cosmos, magnetic fields are present everywhere. For the treatment of fluids and gases coupled to an electromagnetic field, magnetohydrodynamics is usually an effective approximation. Relativistic magnetohydrodynamics is needed to analyze in detail to understand several astrophysical phenomena, including accretion disks, plasma winds, magnetospheres, and astrophysical jets near compact objects where the relativistic effects are considerable. This motivated us to discuss a sonic-supersonic degenerate boundary value problem for a relativistic magnetohydrodynamics system in this article, where the boundary data is prescribed on a characteristic curve and a sonic boundary. We choose a general convex pressure which makes the problem under consideration even more complicated and extensive.

Let us consider the motion of a relativistic fluid in an electromagnetic field in the Minkowski space-time governed by relativistic magnetohydrodynamics equations. In the covariant form, these equations can be written as \cite{anile1987mathematical}
\begin{align}\label{eq: main}
    \begin{cases}
        \partial_j(nu^j)=0, ~~~~~~~~~~~~~~~~~~~~~~~~~~~~~\mathrm{(Conservation ~of ~mass),}\\    \partial_j(\mathbf{T}^{jk})=0,~~~~~~~~~~~~~~~~~~~~~~~~~~~~~\mathrm{(Conservation ~of ~energy-momentum),}\\
        \partial_j(\mathbf{\Psi}^{jk})=0,~~j,k=0, 1, 2, 3~~~\mathrm{(Maxwell ~equations),}\\
    \end{cases}
\end{align}
where the Einstein summation convention has been used, $n$ denotes the average rest-mass density of baryon particles, $u^j$ denotes the four-velocity vector and $\partial_j$ denote the covariant derivative operator with respect to space-time coordinates $(t, x^1, x^2, x^3)$. Further, the tensor $\mathbf{\Psi}^{jk}$ is given by $\mathbf{\Psi}^{jk}=u^jh^k-u^kh^j$ and the energy-momentum tensor $\mathbf{T}^{jk}$ is decomposed into the fluid part $\mathbf{T}_f^{jk}$ and electromagnetic part $\mathbf{T}_m^{jk}$ such that $\mathbf{T}^{jk}= \mathbf{T}_f^{jk}+\mathbf{T}_m^{jk}$ with
\begin{align*}
\begin{cases}
\mathbf{T}_f^{jk}=(p+\rho)u^ju^k+pg^{jk},\\
\mathbf{T}_m^{jk}=|h|^2(u^ju^k+g^{jk}/2)-h^jh^k,
\end{cases}
\end{align*}
where $i=p+\rho$ represents the enthalpy per unit volume such that $p$ is the gas pressure, $\rho=mn+e$ defines the total mass-energy density with $e$ being the specific internal energy and $m$ being the average rest mass. Further, we take $g^{jk}$ as the Minkowski tensor, i.e. $g^{jk}=\pm \mathrm{diag}\{-1, 1, 1, 1\}$ and normalize the speed of light to be one. Also, relations between the four vectors $u^j$ and $h^j$ and the spatial components of the velocity $\Vec{\mathbf{U}}=(u_1, u_2, u_3)$ and the laboratory magnetic field $\Vec{\mathbf{H}}=(H_1, H_2, H_3)$ are $u^j=\gamma(1, \Vec{\mathbf{U}}),~~~ h^j=\gamma\left(\Vec{\mathbf{U}}.\Vec{\mathbf{H}}, \dfrac{\Vec{\mathbf{H}}}{\gamma^2}+\Vec{\mathbf{U}}(\Vec{\mathbf{U}}.\Vec{\mathbf{H}})\right),$ where $\gamma=\dfrac{1}{\sqrt{1-u_1^2-u_2^2-u_3^2}}~(0<u_1^2+u_2^2+u_3^2<1)$ is the Lorentz factor. 
For simplicity of notations and computations to be followed, we rewrite \eqref{eq: main} into the conservative form given by 
\begin{equation}\label{eq: 1.a}
\hspace{-0.35 cm}\begin{cases}
\begin{aligned}
     &(n\gamma)_t+\mathrm{div}~(n\gamma \Vec{\mathbf{U}})=0,\\
     &\hspace{-0.2 cm}\left((i\gamma^2+|\Vec{\mathbf{H}}|^2)\Vec{\mathbf{U}}-(\Vec{\mathbf{U}}.\Vec{\mathbf{H}})\Vec{\mathbf{H}}\right)_t+\mathrm{div}\left(\left((i\gamma^2+|\Vec{\mathbf{H}}|^2)\Vec{\mathbf{U}}-(\Vec{\mathbf{U}}.\Vec{\mathbf{H}})\Vec{\mathbf{H}}\right)\otimes \Vec{\mathbf{U}}-\Vec{\mathbf{H}}\left(\dfrac{\Vec{\mathbf{H}}}{\gamma^2}+(\Vec{\mathbf{U}}.\Vec{\mathbf{H}})\Vec{\mathbf{U}}\right)+p^{tot}\mathbf{I}\right)=0,\\
     &\hspace{-0.2 cm}\left(i\gamma^2-p^{tot}+|\Vec{\mathbf{H}}|^2\right)_t+\mathrm{div}\left((i\gamma^2+|\Vec{\mathbf{H}}|^2)\Vec{\mathbf{U}}-(\Vec{\mathbf{U}}.\Vec{\mathbf{H}})\Vec{\mathbf{H}}\right)=0,\\
     &\Vec{\mathbf{H}}_t-\mathrm{rot}~(\Vec{\mathbf{U}}\times \Vec{\mathbf{H}})=0,\\
     &\mathrm{div} ~\Vec{\mathbf{H}}=0,
\end{aligned}
\end{cases}
\end{equation}
where $\mathbf{I}$ is a $3\times 3$ identity matrix and $p^{tot}=p+p_m$ is the total pressure containing the gas pressure $p$ and the magnetic pressure $p_m$ given by $p_m= \dfrac{1}{2}\left(\dfrac{|\Vec{\mathbf{H}}|^2}{\gamma^2}+(\Vec{\mathbf{U}}.\Vec{\mathbf{H}})^2\right)$. 

Note that $\Vec{\mathbf{U}}\times \Vec{\mathbf{H}}=(u_2 H_3-u_3 H_2, u_3 H_1-u_1H_3, u_1 H_2-u_2H_1)$, which implies that the magnetic field equations can be rewritten as 
\begin{align*}
    \begin{cases}
(H_1)_t-(u_1 H_2-u_2 H_1)_y-(u_3 H_1-u_1 H_3)_z=0,\\
(H_2)_t-(u_2 H_3-u_3 H_2)_z-(u_1 H_2-u_2 H_1)_x=0,\\
(H_3)_t-(u_3 H_1-u_1 H_3)_x-(u_2 H_3-u_3H_2)_y=0.
    \end{cases}
\end{align*}
Therefore, if we assume that the magnetic field is orthogonal to the velocity field such that $\Vec{\mathbf{H}}=(0, 0, H)$ and $\Vec{\mathbf{U}}=(u, v, 0)$ satisfying $\Vec{\mathbf{U}}.\Vec{\mathbf{H}}=0$ and the variables $(u, v, p, \rho, H)$ are independent of the space variable $z$, then the system \eqref{eq: 1.a} is reduced to
\begin{align}\label{eq: 1.2}
    \begin{cases}
(n \gamma)_t+(n \gamma u)_x+(n \gamma v)_y=0,\\
\left((i\gamma^2+H^2)u\right)_t+\left((i\gamma^2+H^2)u^2+p+\dfrac{H^2}{2\gamma^2}\right)_x+\left((i\gamma^2+H^2)uv\right)_y=0,\\
\left((i\gamma^2+H^2)v\right)_t+\left((i\gamma^2+H^2)uv\right)_x+\left((i\gamma^2+H^2)v^2+p+\dfrac{H^2}{2\gamma^2}\right)_y=0,\\
\left(i\gamma^2-p-\dfrac{H^2}{2\gamma^2}+H^2\right)_t+\left((i\gamma^2+H^2)u\right)_x+\left((i\gamma^2+H^2)v\right)_y=0,\\
H_t+(Hu)_x+(Hv)_y=0.
    \end{cases}
\end{align}
From the first and last equation, it is trivial to see that
\begin{align}\label{eq: 1.3}
    \left(\dfrac{H}{n\gamma}\right)_t+u\left(\dfrac{H}{n\gamma}\right)_x+v\left(\dfrac{H}{n\gamma}\right)_y=0,
\end{align}
which implies that $\dfrac{H}{n\gamma}$ is a constant along each stream line. The relation \eqref{eq: 1.3} is essentially an expression of the frozen-in law. To put it another way, all magnetic field lines move with the particles on this line and it appears like magnetic field lines have “frozen” on the fluid moving with them, see \cite{chen2020rarefaction}. For simplicity, we assume that $\dfrac{H}{n\gamma}=\kappa_0$ where $\kappa_0$ is a constant. 

Also, the Rankine-Hugoniot conditions for system \eqref{eq: 1.2} are given by
\begin{align}\label{eq: 1.2a}
    \begin{cases}
        \alpha_t[n\gamma]=\alpha_x[n\gamma u]+\alpha_y[n\gamma v],\\      \alpha_t[(i\gamma^2+H^2)u]=\alpha_x\bigg[(i\gamma^2+H^2)u^2+p+\dfrac{H^2}{2\gamma^2}\bigg]+\alpha_y[(i\gamma^2+H^2)uv],\\
\alpha_t[(i\gamma^2+H^2)v]=\alpha_x[(i\gamma^2+H^2)uv]+\alpha_y\bigg[(i\gamma^2+H^2)v^2+p+\dfrac{H^2}{2\gamma^2}\bigg],\\
\alpha_t\bigg[i\gamma^2-p-\dfrac{H^2}{2\gamma^2}+H^2\bigg]=\alpha_x[(i\gamma^2+H^2)u]+\alpha_y[(i\gamma^2+H^2)v],\\
\alpha_t[H]=\alpha_x[H u]+\alpha_y[H v],\\   
    \end{cases}
\end{align}
where $[(.)]$ denotes the jump in the variable $(.)$ across a discontinuity surface and $(-\alpha_t, \alpha_x, \alpha_y)$ denotes the normal to the discontinuity surface. 

Now one can use the first and fifth equations of \eqref{eq: 1.2a} to obtain 
\begin{align}\label{eq: 1.5a}
    \bigg[\dfrac{H}{n\gamma}\bigg]=0,
\end{align}
which means that $\dfrac{H}{n\gamma}$ remains unchanged across a shock wave. Therefore, in view of \eqref{eq: 1.3} and \eqref{eq: 1.5a}, we consider the following reduced system
\begin{align*}
    \begin{cases}
(n \gamma)_t+(n \gamma u)_x+(n \gamma v)_y=0,\\
\bigg[((i+\kappa_0^2 n^2)\gamma^2)u\bigg]_t+\left((i+\kappa_0^2 n^2)\gamma^2)u^2+p+\dfrac{1}{2}\kappa_0^2 n^2\right)_x+\bigg[(i+\kappa_0^2 n^2)\gamma^2)uv\bigg]_y=0,\\
\bigg[((i+\kappa_0^2 n^2)\gamma^2)v\bigg]_t+\bigg[(i+\kappa_0^2 n^2)\gamma^2)uv\bigg]_x+\left((i+\kappa_0^2 n^2)\gamma^2)v^2+p+\dfrac{1}{2}\kappa_0^2 n^2\right)_y=0,\\
\bigg[((i+\kappa_0^2 n^2)\gamma^2-p-\dfrac{1}{2}\kappa_0^2 n^2)v\bigg]_t+\bigg[(i+\kappa_0^2 n^2)\gamma^2u\bigg]_x+\bigg[(i+\kappa_0^2 n^2)\gamma^2 v\bigg]_y=0.
    \end{cases}
\end{align*}
In this article, we consider that the flow is isentropic and steady which leads to the following system 
\begin{align}\label{eq: 1.1}
\begin{cases}
(n \gamma u)_x+(n \gamma v)_y=0,\\
\bigg[(\rho+p+\kappa_0^2 n^2)\gamma^2)u^2+p^{tot}\bigg]_x+\bigg[((\rho+p+\kappa_0^2 n^2)\gamma^2)uv\bigg]_y=0,\\
\bigg[(\rho+p+\kappa_0^2 n^2)\gamma^2)uv\bigg]_x+\bigg[(\rho+p+\kappa_0^2 n^2)\gamma^2)v^2+p^{tot}\bigg]_y=0,
\end{cases}     
\end{align}
where $0<q=\sqrt{u^2+v^2}<1$ is the flow velocity, $\gamma=(1-q^2)^{-1/2}$ is the normalized Lorentz factor, and $p^{tot}=p+\dfrac{1}{2}\kappa_0^2 n^2$ is the total pressure. Moreover, we suppose that the density $\rho$ and gas pressure $p=p(\rho)$ satisfy the following \cite{chen2018boundary, chen2004stability}
\begin{align}\label{1.5}
    0<\rho<\rho_{\max}<\infty,~~ 0<p(\rho), ~~0<p'(\rho)+\dfrac{\kappa_0^2 n^2}{p+\rho}(\rho)<1, ~~p''(\rho)>0.
\end{align}
Additionally, we assume that for any value of $\rho$, $p''(\rho)$ does not become unbounded, which is true in general for all equations of states of physical importance.

In view of Figure \ref{fig: 1}, we describe the degenerate Cauchy-Goursat boundary value problem for relativistic magnetohydrodynamics precisely as follows:\\\\
\textbf{A Cauchy-Goursat mixed-type boundary value problem:\\\\}
\textit{If $\widehat{BA}$ and $\widehat{BC}$ are two given smooth boundaries of the degenerate mixed-type boundary value problem such that $\widehat{BC}$ is a sonic curve and $\widehat{BA}$ is a positive characteristic curve of the flow. Then develop a local classical supersonic solution to the relativistic magnetohydrodynamics system in the angular domain $ABC$ near the point $B$.}\\

Such a degenerate Cauchy-Goursat problem for two-dimensional steady Euler equations was studied by Li and Hu \cite{li2019degenerate}. They proved the convergence of a sequence of iterations generated by a system of integral equations using the ideas of characteristic decomposition to establish the existence and uniqueness of the solution. Their work was extended to the relativistic Euler equations for polytopic gas by Fan et al. \cite{fan2022sonic} recently. We try to extend these results to a more complicated relativistic magnetohydrodynamics system with a general convex pressure. Note that due to the presence of a general pressure law and magnetic effects, the process of proving the convergence is not straightforward and brings more complexities. The approach to prove the existence and uniqueness of the solution in this article is inspired by the work of Li and Hu \cite{li2019degenerate}. We first establish the existence and uniqueness of the solution in a partial hodograph region by proving the uniform convergence of an iterated sequence generated by a system of integral equations and then obtain the solution to the original problem in the physical plane by transforming the solution obtained in partial hodograph plane via an inverse transformation.

The rest of the article can be arranged as follows. We analyze the basic aspects of the steady isentropic rotational relativistic magnetohydrodynamic equations \eqref{eq: 1.1} and derive the characteristic decompositions of angle variables in section 2. In section 3, we precisely define the main problem and specify the boundary data on the characteristic curve and sonic curve, respectively. The purpose of section 4 is to establish the existence and uniqueness of solutions in the partial hodograph plane. We use an inverse transformation to convert the solution obtained in the partial hodograph plane back into the physical plane and verify that they satisfy our main problem in section 5.

}
\section{Basic aspects of steady isentropic irrotational  relativistic magnetohydrodynamics}\label{2}
Assuming that the relativistic flow is irrotational, i.e., $u_y=v_x$ and using first equation of \eqref{eq: 1.1} and $\gamma_x=\gamma^3(uu_x+vv_x), ~~\gamma_y=\gamma^3(uu_y+vv_y),$ we can write the second equation of system \eqref{eq: 1.1} as follows
\begin{equation}\label{eq: 2.1}
    \dfrac{n}{\gamma}\bigg\{\dfrac{i+\frac{H^2}{\gamma^2}}{n}\gamma_x+\dfrac{\gamma}{n}p_x^{tot}+\dfrac{\gamma u}{n}\bigg[n\gamma u\left(\dfrac{\gamma\left( i+\frac{H^2}{\gamma^2}\right)}{n}\right)_x+n\gamma v\left(\dfrac{\gamma \left(i+\frac{H^2}{\gamma^2}\right)}{n}\right)_y\bigg]\bigg\}=0.
\end{equation}
Further, by the first law of thermodynamics, we have
\begin{align}\label{eq: 2.2}
    d\left(\dfrac{i}{n}\right)=\dfrac{1}{n}dp+TdS,
\end{align}
where $S$ is the entropy of the flow and $T$ denotes absolute temperature. Since the entropy is assumed to be constant, in view of \eqref{eq: 2.2}, we have
\begin{align}\label{eq: 2.3}
    d\left(\dfrac{i}{n}\right)=\dfrac{1}{n}dp
\end{align}
or equivalently
\begin{align*}
    \dfrac{dn}{n}=\dfrac{d\rho}{p+\rho},
\end{align*}
i.e., 
\begin{align}
    n=n(\rho)=n_0 e^{\displaystyle \int_{\rho_0}^{\rho} \dfrac{ds}{s+p(s)}},
\end{align}
where $n_0$ is the proper number density of baryon particles when the mass-energy density of the fluid is $\rho_0>0$. It is straightforward to observe that 
\begin{align}
\begin{cases}
    n'(\rho)=\dfrac{n(\rho)}{p(\rho)+\rho}>0,\\
    n''(\rho)=-\dfrac{p'(\rho)n(\rho)}{(p(\rho)+\rho)^2}<0.
\end{cases}
\end{align}

Also, noting that $H=\kappa_0 n \gamma$, one has
\begin{align}\label{eq: 2.6}
    d\left(\dfrac{H^2}{n\gamma^2}\right)=\dfrac{dp_m}{n},
\end{align}
where $p_m$ is the magnetic pressure defined by $p_m=\dfrac{H^2}{2\gamma^2}$. Clearly, one can use \eqref{eq: 2.3} and \eqref{eq: 2.6} to obtain
\begin{align}\label{eq: 2.7}
    d\left(\dfrac{i+\frac{H^2}{\gamma^2}}{n}\right)=\dfrac{dp^{tot}}{n}.
\end{align}
Thus, exploiting \eqref{eq: 2.1} we have
\begin{align}\label{eq: a}
    \left(\dfrac{\gamma \left(i+\frac{H^2}{\gamma^2}\right)}{n}\right)_x+\dfrac{\gamma u}{n}\bigg[n\gamma u\left(\dfrac{\gamma \left(i+\frac{H^2}{\gamma^2}\right)}{n}\right)_x+n\gamma v\left(\dfrac{\gamma \left(i+\frac{H^2}{\gamma^2}\right)}{n}\right)_y\bigg]=0.
\end{align}
In a similar fashion, using the last equation of system \eqref{eq: 1.1}, one can get 
\begin{align}\label{eq: b}
    \left(\dfrac{\gamma \left(i+\frac{H^2}{\gamma^2}\right)}{n}\right)_y+\dfrac{\gamma v}{n}\bigg[n\gamma u\left(\dfrac{\gamma \left(i+\frac{H^2}{\gamma^2}\right)}{n}\right)_x+n\gamma v\left(\dfrac{\gamma \left(i+\frac{H^2}{\gamma^2}\right)}{n}\right)_y\bigg]=0.
\end{align}
Now one can easily observe that the matrix $\begin{bmatrix}
    1+\gamma^2u^2&\gamma^2 uv\\
    \gamma^2uv&1+\gamma^2v^2
    \end{bmatrix}$ is nonsingular which implies that $\left(\dfrac{\gamma \left(i+\frac{H^2}{\gamma^2}\right)}{n}\right)_x=\left(\dfrac{\gamma \left(i+\frac{H^2}{\gamma^2}\right)}{n}\right)_y=0$ and yields the Bernoulli's law for two-dimensional steady relativistic magnetohydrodynamics of the form
\begin{align}\label{eq: 2.10}
\dfrac{\gamma \left(i+\frac{H^2}{\gamma^2}\right)}{n}= \mathrm{constant.}
\end{align}
Without loss of generality, we rewrite \eqref{eq: 2.10} as:
\begin{align}\label{eq: 2.11}
\dfrac{\gamma \left(i+\frac{H^2}{\gamma^2}\right)}{n}= m\hat{\gamma},
\end{align}
where $\hat{\gamma}^{-1}=\sqrt{1-\hat{q}^2}$ is a constant.
\begin{l1}
If for all $n>0$, $\dfrac{\partial p}{\partial n}>0$  then there exists a constant $\hat{q}~(0<\hat{q}<1)$ usually referred as the limit speed of the flow such that $q<\hat{q}$, where $q$ is the flow speed. Further, $q$ approaches $\hat{q}$ when $n$ approaches $0$.
\end{l1}
\begin{proof}
We use \eqref{eq: 2.7} to obtain
$$\dfrac{d}{dn}\left(\dfrac{i+\frac{H^2}{\gamma^2}}{n}\right)=\dfrac{d\left(\dfrac{i+\frac{H^2}{\gamma^2}}{n}\right)}{dp^{tot}}.\dfrac{dp^{tot}}{dn}=\dfrac{1}{n}.\dfrac{dp^{tot}}{dn}>0$$ 
for $n>0$.

In view of \eqref{eq: 2.11} and the fact that $\dfrac{i+\dfrac{H^2}{\gamma^2}}{n}=\dfrac{mn+e+p+\kappa_0^2 n^2}{n}\geq m$ for $n\geq 0$, it is straightforward to observe that $q<\hat{q}<1$ and $q\longrightarrow \hat{q}$ as $n\longrightarrow 0$.
\end{proof}
Again by employing \eqref{eq: 2.11}, system \eqref{eq: 1.1} can be converted into
\begin{align}\label{eq: 2.10 a}
\begin{cases}
    \gamma\Bigg[\left(\left(i+\frac{H^2}{\gamma^2}\right)\gamma u\right)_x+\left(\left(i+\frac{H^2}{\gamma^2}\right)\gamma v\right)_y\Bigg]+\left(i+\frac{H^2}{\gamma^2}\right)\gamma [u\gamma_x+v\gamma_y]=0,\\
    \gamma u\Bigg[\left(\left(i+\frac{H^2}{\gamma^2}\right)\gamma u\right)_x+
    \left(\left(i+\frac{H^2}{\gamma^2}\right)\gamma v\right)_y\Bigg]+\left(i+\frac{H^2}{\gamma^2}\right)\gamma [u(\gamma u)_x+v(\gamma u)_y]+p^{tot}_x=0,\\
    \gamma v\Bigg[\left(\left(i+\frac{H^2}{\gamma^2}\right)\gamma u\right)_x+\left(\left(i+\frac{H^2}{\gamma^2}\right)\gamma v\right)_y\Bigg]+\left(i+\frac{H^2}{\gamma^2}\right)\gamma [u(\gamma v)_x+v(\gamma v)_y]+p^{tot}_y=0.
\end{cases}    
\end{align}
Further, by considering scalar product of $(\gamma, -\gamma u, -\gamma v)$ with \eqref{eq: 2.10 a} leads to
\begin{align}\label{eq: 2.13}
    \left(i+\frac{H^2}{\gamma^2}\right)(\gamma u)_x+\left(i+\frac{H^2}{\gamma^2}\right) (\gamma v)_y=-\gamma u \rho_x-\gamma v \rho_y.
\end{align}
Moreover, from the second and third equation of \eqref{eq: 1.1}, one can obtain
\begin{align}\label{eq: 2.14}
    \begin{cases}
        \left(i+\frac{H^2}{\gamma^2}\right)\gamma^2 u u_x+\left(i+\frac{H^2}{\gamma^2}\right) \gamma^2 vu_y+p_x^{tot}=0,\\
        \left(i+\frac{H^2}{\gamma^2}\right)\gamma^2 u v_x+\left(i+\frac{H^2}{\gamma^2}\right) \gamma^2 vv_y+p_y^{tot}=0.
    \end{cases}
\end{align}
Again by taking the scalar product of $(\gamma u, \gamma v)$ with \eqref{eq: 2.14}, one can obtain
\begin{align}\label{eq: 2.15}
    \left(i+\frac{H^2}{\gamma^2}\right)\gamma^3[u^2 u_x+uv(u_y+v_x)+v^2v_y]=-w^2(\gamma u \rho_x+\gamma v \rho_y),
\end{align}
where $w=\sqrt{p'(\rho)+\dfrac{\kappa_0^2 n^2}{p+\rho}}$ is the magneto-acoustic sound.

Therefore, we combine \eqref{eq: 2.13} and \eqref{eq: 2.15} to convert the governing system \eqref{eq: 1.1} in the form:
\begin{align}\label{eq: 2.16}
    \begin{cases}
    (M_1^2-1)u_x+M_1M_2(v_x+u_y)+(M_2^2-1)v_y=0,\\
    u_y-v_x=0,
    \end{cases}
\end{align}
where $M_1=\dfrac{\gamma u}{w\gamma_w}, ~~M_2=\dfrac{\gamma v}{w\gamma_w},$ and $\gamma_w=\dfrac{1}{\sqrt{1-w^2}}$.
 
For smooth solutions, \eqref{eq: 2.16} can be rewritten in the matrix form as
\begin{align}\label{eq: 2.17}
    \begin{bmatrix}
    M_1^2-1&M_1M_2\\
    0&-1
    \end{bmatrix}
    \begin{bmatrix}
    u\\
    v
    \end{bmatrix}_x
    +
    \begin{bmatrix}
    M_1M_2&M_2^2-1\\
    1&0
    \end{bmatrix}
    \begin{bmatrix}
    u\\
    v
    \end{bmatrix}_y=
    \begin{bmatrix}
    0\\0
    \end{bmatrix}.
\end{align}
The eigenvalues of the system \eqref{eq: 2.17} are  $\Lambda_\pm=\dfrac{M_1M_2\pm \sqrt{M_1^2+M_2^2-1}}{M_1^2-1}$ and the left eigenvectors are $l_\pm=(1, \mp \sqrt{M^2-1})$, such that $M=\sqrt{M_1^2+M_2^2}=\dfrac{\gamma q}{\gamma_w w}$ is  the proper Mach number. Clearly, \eqref{eq: 2.17} is a mixed-type system and may change its behavior across the sonic curve ($M=1$). The system \eqref{eq: 2.17} is hyperbolic for $M>1$ while for $M<1$ system \eqref{eq: 2.17} changes its nature to elliptic. Moreover, one can define the families of wave characteristics as follows
\begin{equation}
    \dfrac{dy}{dx}=\Lambda_\pm.
\end{equation}
Furthermore, one can easily obtain the characteristic equations of the form 
\begin{align}\label{eq: 2.18}
\begin{cases}
{\partial}_{+}u+\Lambda_- {\partial}_{+}v=0,\\
{\partial}_{-}u+\Lambda_+ {\partial}_{-}v=0,
\end{cases}
\end{align} 
where ${\partial}_{\pm}=\partial_x+\Lambda_{\pm}\partial_y$. 
\subsection{Characteristic decompositions of the flow variables}
Differentiating Bernoulli's law \eqref{eq: 2.11} with respect to $q$ gives
\begin{align}\label{eq: 2.20}
    \dfrac{\gamma^3 q \left(i+\frac{H^2}{\gamma^2}\right)}{n}+\gamma\dfrac{d(\left(i+\frac{H^2}{\gamma^2}\right)/n)}{dp^{tot}}\dfrac{dp^{tot}}{dw}\dfrac{dw}{dq}=0.
\end{align}
Then by noting that $w^2=p'(\rho)+\dfrac{\kappa_0^2 n^2}{p+\rho}$ we have $\dfrac{dp^{tot}}{dw}=\dfrac{2i^2w^3}{i^2p''(\rho)+\kappa_0^2 n^2(1-p'(\rho))}$, which is used in \eqref{eq: 2.20} along with \eqref{eq: 2.7} to yield
\begin{align}\label{eq: 2.20a}
    \dfrac{dw}{dq}=-\dfrac{\left(i+\frac{H^2}{\gamma^2}\right)q\gamma^2\bigg[ i^2 p''(\rho)+\kappa_0^2 n^2(1-p'(\rho))\bigg]}{2i^2 w^3}<0.
\end{align}
Also, by $M=\dfrac{\gamma q}{w\gamma_w}$, it is easy to see that
\begin{align*}
    \dfrac{dM}{dq}=M\left(\dfrac{1}{q(1-q^2)}-\dfrac{1}{w(1-w^2)}\dfrac{dw}{dq}\right).
\end{align*}
Then noting that $\dfrac{dw}{dq}<0$ and $0<w^2, q^2<1$ we have $\dfrac{dM}{dq}>0$.

Furthermore, we employ the following expressions of velocity \cite{chen2018boundary} 
 \begin{equation}\label{eq: 2.21}
 u=\dfrac{w\gamma_w}{\gamma}\dfrac{\cos\theta}{\sin \omega}, ~~v=\dfrac{w\gamma_w}{\gamma}\dfrac{\sin\theta}{\sin \omega},
 \end{equation}
and employ the normalized directional derivatives of the form  
\begin{align}\label{eq: 2.22}
\begin{cases}
    &\bar{\partial}_+= \cos \alpha \dfrac{\partial}{\partial x}+ \sin \alpha \dfrac{\partial}{\partial y}, ~~\bar{\partial}_-= \cos \beta \dfrac{\partial}{\partial x}+ \sin \beta \dfrac{\partial}{\partial y},\vspace{0.2 cm}\\
    &\bar{\partial}_0=\cos \theta \dfrac{\partial}{\partial x}+ \sin \theta \dfrac{\partial}{\partial y},~~\bar{\partial}_\perp= -\sin \theta \dfrac{\partial}{\partial x}+ \cos \theta \dfrac{\partial}{\partial y}.
    \end{cases}
\end{align}
From \eqref{eq: 2.22}, we have
\begin{align}
\dfrac{\partial}{\partial x}=-\dfrac{\sin \beta \bar{\partial}_+-\sin \alpha \bar{\partial}_-}{ \sin (2\omega)},~~\dfrac{\partial}{\partial y}=\dfrac{\cos \beta \bar{\partial}_+-\cos \alpha \bar{\partial}_-}{ \sin (2\omega)},~~\bar{\partial}_0=\dfrac{\bar{\partial}_++\bar{\partial}_-}{2\cos \omega},~~\bar{\partial}_\perp=\dfrac{\bar{\partial}_+-\bar{\partial}_-}{2\sin \omega}.
\end{align}
Therefore, a direct computation leads to the following decompositions of $u$ and $v$
 \begin{align}\label{eq: 2.24}
     \begin{cases}
     \bar{\partial}_\pm u=\dfrac{\gamma_w}{\gamma \sin^2 \omega}\bigg[ f(w) \sin \omega \cos \theta \bar{\partial}_\pm w-w \sin \omega \sin \theta \bar{\partial}_\pm \theta-w  \cos \theta \cos \omega \bar{\partial}_\pm\omega\bigg],\vspace{0.2 cm}\\
    \bar{\partial}_\pm v=\dfrac{\gamma_w}{\gamma \sin^2 \omega}\bigg[  f(w)\sin \omega \sin \theta  \bar{\partial}_\pm w+w   \sin \omega \cos \theta \bar{\partial}_\pm \theta-w  \sin \theta \cos \omega \bar{\partial}_\pm\omega\bigg],
     \end{cases}
 \end{align}
 where $0<f(w)=\left(\gamma_w^2+\dfrac{2i^2 w^4}{\left(i+\frac{H^2}{\gamma^2}\right)\bigg[ i^2p''(\rho)+\kappa_0^2 n^2(1-p'(\rho))\bigg]}\right)<\infty$.
 
Utilizing \eqref{eq: 2.24} and performing direct calculations leads us to
 \begin{align}
 \begin{cases}
     \bar{\partial}_+ \omega&=\dfrac{wF_1(w, \omega)}{2\left(i+\frac{H^2}{\gamma^2}\right)q\gamma \bigg[ i^2 p''(\rho)+\kappa_0^2 n^2(1-p'(\rho))\bigg]\cos \omega \gamma_w}\dfrac{\bar{\partial}_+w}{w},\\
     \bar{\partial}_- \omega&=\dfrac{wF_1(w, \omega)}{2\left(i+\frac{H^2}{\gamma^2}\right)q\gamma\bigg[ i^2 p''(\rho)+\kappa_0^2 n^2(1-p'(\rho))\bigg]\cos \omega \gamma_w}\dfrac{\bar{\partial}_-w}{w},\\
     \bar{\partial}_+ \theta&=\dfrac{(F_2(w, \omega)-wF_1(w, \omega)\cos^2\omega)}{2\left(i+\frac{H^2}{\gamma^2}\right)q\gamma \bigg[ i^2 p''(\rho)+\kappa_0^2 n^2(1-p'(\rho))\bigg]\cos \omega\sin^2\omega \gamma_w}\dfrac{\bar{\partial}_+w}{w},\\
      \bar{\partial}_- \theta&=\dfrac{(wF_1(w, \omega)\cos^2\omega-F_2(w, \omega))}{2\left(i+\frac{H^2}{\gamma^2}\right)q\gamma \bigg[ i^2 p''(\rho)+\kappa_0^2 n^2(1-p'(\rho))\bigg]\cos \omega\sin^2\omega \gamma_w}\dfrac{\bar{\partial}_-w}{w},
      \end{cases}
 \end{align}
 where 
 \begin{align}
 \begin{cases}
 F_1(a, \omega)=2\left(i+\frac{H^2}{\gamma^2}\right)\bigg[i^2 p''(\rho)+\kappa_0^2 n^2(1-p'(\rho)) \bigg]\gamma_w^2f(w)+4 i^2 w^2\sin^2\omega> 0,\\
 F_2(a, \omega)=\left(i+\frac{H^2}{\gamma^2}\right)q\gamma\gamma_w f(w)\bigg[i^2 p''(\rho)+\kappa_0^2 n^2(1-p'(\rho))\bigg]\cos \omega \sin 2\omega.
 \end{cases}
 \end{align}
Noting the expressions of velocity components from \eqref{eq: 2.21} and employing \eqref{eq: 2.24} and \eqref{eq: 2.20a}, it is easy to see that $\omega=\omega(\rho)$. Therefore, one can use $\dfrac{d\omega}{dw}=\dfrac{d\omega}{dM}\dfrac{dM}{dq}\dfrac{dq}{dw}\dfrac{dw}{d\rho}$ and $\dfrac{dw}{dq}<0,~\dfrac{dM}{dq}>0$ along with $\sin \omega=\dfrac{1}{M}$ to obtain $\dfrac{d\omega}{dw}>0$. Further, we denote $\sin \omega=\dfrac{1}{M}:=\varpi, \omega\in [\varpi_0, \pi/2],$ ($\varpi_0>0$ is a constant). Thus by utilizing the inverse function theorem, one can observe that $\rho=\rho(\omega)=\rho(\sin^{-1}\varpi)$ or in other words $w=w(\omega)=w(\sin^{-1}\varpi)$. Thus we denote $F_1(w, \omega)=:F_1(\varpi)$ and $F_2(w, \omega)=:F_2(\varpi)$ to get
 \begin{align}\label{eq: 2.27}
     \begin{cases}
     \bar{\partial}_+\theta+\dfrac{4i^2(\varpi)w^2(\varpi)\cos \omega}{F_1(\varpi)}\bar{\partial}_+\varpi=0,\\
     \bar{\partial}_-\theta-\dfrac{4i^2(\varpi)w^2(\varpi) \cos \omega}{F_1(\varpi)}\bar{\partial}_-\varpi=0.
     \end{cases}
 \end{align}
In view of the continuity of $\dfrac{2i^2(s)w^2(s)}{sF_1(s)}, s\in [\sin \varpi_0, 1]$, we set
 \begin{align*}
     I(\varpi)=\int_{\sin k_0}^{\varpi}\dfrac{2i^2(s)w^2(s)}{sF_1(s)}ds, ~~\varpi\in [\sin \varpi_0, 1].
 \end{align*}
Thus, we denote $I:=I(\varpi)$ to reduce \eqref{eq: 2.27} as
 \begin{align}\label{eq: 2.28}
     \begin{cases}
     \bar{\partial}_+\theta+\sin 2\omega \bar{\partial}_+I=0,\\
     \bar{\partial}_-\theta-\sin 2\omega \bar{\partial}_-I=0
     \end{cases}
 \end{align}
 with 
 \begin{align}\label{eq: 2.29}
     \bar{\partial}_i \varpi=\dfrac{\varpi F_1(\varpi)}{2i^2(\varpi)w^2(\varpi)}\bar{\partial}_i I, ~~i=0, \pm.
 \end{align}
 Moreover, we make use of the commutator relation given by \cite{li2019degenerate}
 \begin{align*}
  \bar{\partial}_-\bar{\partial}_+ -\bar{\partial}_+\bar{\partial}_-&=\dfrac{1}{\sin 2\omega}\bigg[(\cos 2 \omega \bar{\partial}_-\alpha-\bar{\partial}_+\beta)\bar{\partial}_+ -(\bar{\partial}_-\alpha-\cos 2\omega \bar{\partial}_+\beta)\bar{\partial}_- \bigg],
 \end{align*}
which provides the second-order decomposition of $I$ and $\theta$ of the form
\begin{align*}
 \begin{cases}
   \bar{\partial}_-\bar{\partial}_+I-\bar{\partial}_+\bar{\partial}_-I&=\dfrac{1}{\sin 2\omega}\bigg[(\cos 2 \omega \bar{\partial}_-\alpha-\bar{\partial}_+\beta)\bar{\partial}_+ I -(\bar{\partial}_-\alpha-\cos 2\omega \bar{\partial}_+\beta)\bar{\partial}_- I\bigg], \vspace{0.2 cm}\\
   \bar{\partial}_-\bar{\partial}_+\theta-\bar{\partial}_+\bar{\partial}_-\theta&=\dfrac{1}{\sin 2\omega}\bigg[(\cos 2 \omega \bar{\partial}_-\alpha-\bar{\partial}_+\beta)\bar{\partial}_+ \theta-(\bar{\partial}_-\alpha-\cos 2\omega \bar{\partial}_+\beta)\bar{\partial}_- \theta \bigg].
    \end{cases} 
 \end{align*}
For subsequent discussions in the article, we define $f(w):=f(w(\varpi))$, $W:=\bar{\partial}_+I$, $Z:=\bar{\partial}_-I$, and employ $F_1(\varpi)=2\left(i+\frac{H^2}{\gamma^2}\right)\left(i^2 p''(\rho)+\kappa_0^2 n^2(1-p'(\rho))\right)\gamma_{w(\varpi)}^2f(w(\varpi))+4i^2(\varpi)w^2(\varpi)\varpi^2$ to derive the following characteristic decompositions of $W$ and $Z$:
 \begin{align}\label{eq: 2.30}
     \begin{cases}
         &\hspace{-0.5 cm}\bar{\partial}_-W=W\Bigg[W+Z\cos 2\omega +\dfrac{F_1(\varpi)}{4i^2(\varpi)w^2(\varpi)\cos^2 \omega}(W-Z\cos 2\omega)\Bigg],\vspace{0.2 cm}\\
    &\hspace{-0.5 cm} \bar{\partial}_+Z=Z\Bigg[Z+W\cos 2\omega +\dfrac{F_1(\varpi)}{4i^2(\varpi)w^2(\varpi)\cos^2 \omega}(Z-W\cos 2\omega)\Bigg],
     \end{cases}
 \end{align}
 or
 \begin{align}\label{eq: 2.31}
     \begin{cases}
         &\hspace{-0.5 cm}\bar{\partial}_-W=\dfrac{W}{\cos^2 \omega}\Bigg[W+Z\cos^2 2\omega +\dfrac{\left(i+\frac{H^2}{\gamma^2}\right)\bigg[ i^2 p''(\rho)+\kappa_0^2 n^2(1-p'(\rho))\bigg]\gamma_{w}^2f(w)}{2i^2(\varpi)w^2(\varpi)}(W-Z\cos 2\omega)\Bigg],\vspace{0.2 cm}\\
    &\hspace{-0.5 cm} \bar{\partial}_+Z=\dfrac{Z}{\cos^2 \omega}\Bigg[Z+W\cos^2 2\omega +\dfrac{\left(i+\frac{H^2}{\gamma^2}\right)\bigg[ i^2 p''(\rho)+\kappa_0^2 n^2(1-p'(\rho))\bigg]\gamma_{w}^2f(w)}{2i^2(\varpi)w^2(\varpi)}(Z-W\cos 2\omega)\Bigg].
     \end{cases}
 \end{align}
\section{The boundary data and main result}
 We now provide the boundary data for our main problem on the boundaries $\widehat{BC}$ and $\widehat{AB}$. Assume that $\widehat{AB}: x=\psi(y), ~y\in [y_A, y_B]$ and $\widehat{BC}: y=\varphi(x), x\in [x_B, x_C]$ are two smooth curves and $(\theta, \varpi)|_{\widehat{BC}}=(\hat{\theta}, \hat{\varpi})(x)$ and $(\theta, \varpi)|_{\widehat{AB}}=(\tilde{\theta}, \tilde{\varpi})(y)$ are given boundary data on the curves $\widehat{BC}$ and $\widehat{AB}$ which satisfy
\begin{align}\label{eq: 3.1}
    \varphi(x)\in C^3([x_B, x_C)), ~~\hat{\theta}(x)\in C^3([x_B, x_C)),~~\hat{\varpi}(x)=1~~ \forall~~ x\in [x_B, x_C]
\end{align}
and 
\begin{align}\label{eq: 3.2}
\begin{cases}
    x_B=\psi(y_B), ~~\psi(y)\in C^4((y_A, y_B]), ~~\tilde{\theta}(y)= \mathrm{\cot^{-1}} \psi'(y)-\sin^{-1} \tilde{\varpi}(y),~~\tilde{\varpi}(y_B)=1,\vspace{0.3 cm}\\
    \dfrac{1}{\sqrt{1-\tilde{\varpi}^2(y_A)}}-\dfrac{4i^2(\tilde{\varpi})w^2(\tilde{\varpi})\sqrt{1-\tilde{\varpi}^2(y_A)}}{F_1(\tilde{\varpi}(y_A))}>0, \vspace{0.3 cm}\\
    ~~\sin^{-1} \tilde{\varpi}(y)-Q(\tilde{\varpi}(y))= \mathrm{\cot^{-1}} \psi'(y)-\big[\hat{\theta}(x_B)+Q(\hat{\varpi}(x_B))\big],
\end{cases}    
\end{align}
where $$Q(\varpi)=\displaystyle \int \dfrac{4i^2(\varpi)w^2(\varpi) \sqrt{1-\varpi^2}}{F_1(\varpi)} d\varpi .$$
The last equality of \eqref{eq: 3.2} indeed holds. From the first equation of \eqref{eq: 2.27} it is clear that $[\tilde{\theta}+Q(\tilde{\varpi})]|_{\widehat{AB}}=\hat{\theta}(x_B)+Q(\hat{\varpi}(x_B))$. Consequently $(\tilde{\theta}, \tilde{\varpi})(y)$ satisfy the last equality of \eqref{eq: 3.2}.
\subsection{Boundary data of $(W, Z)$}
Now we proceed to provide the boundary data of $W$ and $Z$ on the curves $\widehat{AB}$ and $\widehat{BC}$ using the boundary information of $(\theta, \varpi)$.

Using \eqref{eq: 3.2}, we obtain
\begin{align*}
    \tilde{\theta}(B)= \mathrm{\cot^{-1}}\psi'(B)-\sin^{-1} \tilde{\varpi}(B)=\hat{\theta}(B)+Q(\hat{\varpi}(B))-Q(\tilde{\varpi}(B))=\hat{\theta}(B),
\end{align*}
which together with $\hat{\varpi}(B)=\tilde{\varpi}(B)=1$ gives $(\hat{\theta}, \hat{\varpi})(x_B)=(\tilde{\theta}, \tilde{\varpi})(y_B)$.  In order to check that the compatibility conditions hold true at point $B$, we need to verify that the equations in \eqref{eq: 2.27} are satisfied at point $B$. Moreover, one can easily check that the second equation of \eqref{eq: 2.27} is satisfied at point $B$. For this, we just need to check that $\bar{\partial}_- \theta|_{B}=0$ as $\cos \omega|_{B}=0$.

Now we use \eqref{eq: 3.2} to obtain
\begin{align}\label{eq: 3.3a}
    \tilde{\theta}'(y)=-\dfrac{4i^2(\tilde{\varpi})w^2(\tilde{\varpi})\tilde{\varpi}'(y)\sqrt{1-\tilde{\varpi}^2(y)}}{F_1(\tilde{\varpi}(y))},
\end{align}
which implies that $\tilde{\theta}'(B)=0$ or in other words $\bar{\partial}_+\theta=0$ at $B$. Further, since $\cos \omega|_{B}=0$, therefore noting the expression of $\bar{\partial}_0 \theta$ we must have $\bar{\partial}_-\theta (B)=-\bar{\partial}_+\theta (B)=0$, which proves that the compatibility condition is satisfied at $B$.

Now we discuss the boundary values of $W$ and $Z$ on the curves $\widehat{AB}$ and $\widehat{BC}$. Noting that $\widehat{AB}: x=\psi(y)$ is a positive characteristic, we must have
\begin{align}\label{eq: 3.3}
    W|_{\widehat{AB}}=\dfrac{2i^2(\tilde{\varpi}(y))w^2(\tilde{\varpi}(y))\tilde{\varpi}'(y)}{\tilde{\varpi}(y)F_1(\tilde{\varpi}(y))\sqrt{1+(\psi'(y))^2}}=: \tilde{b}_0(y)\in C^2((y_A, y_B]). 
\end{align}
Again since $\cos \omega|_{\widehat{BC}}=0$, we must have $W+Z|_{\widehat{BC}}=2\cos\omega \bar{\partial}_0 I=0$.  Also, by exploiting \eqref{eq: 2.28} we obtain
\begin{align*}
    W-Z=-\dfrac{\bar{\partial}_+\theta+\bar{\partial}_-\theta}{\sin (2\omega)}=-\dfrac{\bar{\partial}_0\theta}{\sin \omega},
\end{align*}
which yields
\begin{align}
    W|_{\widehat{BC}}=-Z|_{\widehat{BC}}=-\dfrac{\bar{\partial}_0 \theta}{2}\bigg|_{\widehat{BC}}.
\end{align}
Therefore, in order to find the boundary values of $W$ and $Z$, we must derive the boundary information of $\bar{\partial}_0 \theta$ on the curve $\widehat{BC}$. Again from \eqref{eq: 2.27}, one can observe that $\bar{\partial}_\perp \theta=-(W+Z)\cos \omega$, which implies $\bar{\partial}_\perp \theta|_{\widehat{BC}}=0$. Therefore, in view of $\theta(x, y)|_{\widehat{BC}}=\theta(x, \varphi(x))=\hat{\theta}(x)$, we must have
\begin{align*}
    \hat{\theta}'(x)=\theta_x(x, \varphi(x))+\varphi'(x)\theta_y(x, \varphi(x)),
\end{align*}
which together with $\bar{\partial}_{\perp}\theta=0$ implies
\begin{align*}
    \theta_x(x, \varphi(x))=\dfrac{\hat{\theta}'\cos \hat{\theta}'}{\cos \hat{\theta}'+\varphi' \sin \hat{\theta}'},~~ \theta_y(x, \varphi(x))=\dfrac{\hat{\theta}'\sin \hat{\theta}'}{\cos \hat{\theta}'+\varphi' \sin \hat{\theta}'},
\end{align*}
or in other words
\begin{align*}
\bar{\partial}_0\theta|_{\widehat{BC}}=\dfrac{\hat{\theta}'}{\cos \hat{\theta}'+\varphi' \sin \hat{\theta}'}.    
\end{align*}
which implies that 
\begin{align}\label{eq: 3.5}
    W|_{\widehat{BC}}=-Z|_{\widehat{BC}}=-\dfrac{\hat{\theta}'}{2(\varphi' \sin \hat{\theta}'+\cos \hat{\theta}')}=: -\hat{a}_0(x)\in C^2([x_B, x_C)). 
\end{align}
Furthermore, the above derivation process clearly shows that $-\hat{a}_0(x_B)=\tilde{b}_0(y_B)$. Considering the continuity of the functions $\hat{a}_0(x_B)$ and $\tilde{b}_0(y_B)$ there must exist two sufficiently small positive constants $\epsilon_0, \delta_0$ such that $\tilde{b}_0(y)\geq \epsilon_0$ for arbitrary $y\in [y_B-\delta_0, y_B]$ and $-\hat{a}_0(x_B)\geq \epsilon_0$ for arbitrary $x\in [x_B, x_B+\delta_0]$. 

One can take into account the following assumptions; given the solutions are being considered in the neighbourhood of point $B$.
\begin{align}\label{eq: 3.6}
    \hat{\theta}'(x)\leq -\epsilon_0,~~\hat{a}_0(x)\leq -\epsilon_0,~~ \forall~ x\in [x_B, x_C), ~~\psi''(y)\leq \epsilon_0,~ \tilde{b}_0(y)\geq \epsilon_0,~~\forall~ y\in (y_A, y_B].
\end{align}
\section{Solution in a partial hodograph plane}
In this section, we develop the existence and uniqueness of a solution to the nonlinear system \eqref{eq: 2.31} with the boundary data \eqref{eq: 3.3} and \eqref{eq: 3.5} under the assumptions \eqref{eq: 3.6} in the neighbourhood of the point $B$ by utilizing a partial hodograph mapping. We convert the system \eqref{eq: 2.31} into a system of linear singular equations by a change of variables and then prove the convergence of a sequence of iterations generated by a set of integral equations.
\subsection{Formulation of the main problem in $(t, r)$ plane}
We now reformulate the main problem by using a change of variables. First we define a partial hodograph mapping $(x, y)\longrightarrow (t, r)$ of the form
\begin{align}\label{eq: 4.1}
    t=\cos \omega(x, y),~~ r=\theta(x, y).
\end{align}
It is a simple observation to see that $\varpi=\sqrt{1-t^2}$. Thus, we denote $F_1(\varpi):=\hat{F}_1(t)$ where
\begin{align}
    \hat{F}_1(t)=4i^2(t)w^2(t)(1-t^2)+2\left(i+\frac{H^2}{\gamma^2}\right)\bigg[i^2(t)p''(\rho(t))+\kappa_0^2 n^2(t)(1-p'(\rho(t))) \bigg]\gamma_w^2(t)f(w(t))>0.
\end{align}
Further, by exploiting \eqref{eq: 2.28} and \eqref{eq: 2.29}, we obtain
\begin{align}\label{eq: 4.3}
    J&:=\dfrac{\partial(t, r)}{\partial(x, y)}=\begin{vmatrix}
    \dfrac{\partial t}{\partial x}& \dfrac{\partial t}{\partial y}\vspace{0.2 cm}\\
    \dfrac{\partial r}{\partial x}& \dfrac{\partial r}{\partial y}
    \end{vmatrix}\nonumber\\
    &=\dfrac{(1-t^2)\hat{F}_1(t)}{i^2(t)w^2(t)t}WZ\neq 0~~\mathrm{for}~0\leq t<1.
\end{align}
Also, an easy computation now yields
\begin{align}\label{eq: 4.4}
    \bar{\partial}_+=-\dfrac{2F(t)}{t}W\partial_t-2\sqrt{1-t^2}Wt\partial_r,~~\bar{\partial}_-=-\dfrac{2F(t)}{t}Z\partial_t+2\sqrt{1-t^2}Zt\partial_r,
\end{align}
where $F(t)=\dfrac{(1-t^2)\hat{F}_1(t)}{4i^2(t)\hat{w}^2(t)}>1.$

As a result, system \eqref{eq: 2.31} can be represented as a semi-linear system in the $(t, r)$ plane, which takes the form
\begin{align}\label{eq: 4.5}
\begin{cases}
    &\hspace{-0.5 cm}W_t-\dfrac{t^2\sqrt{1-t^2}}{F(t)}W_r=-\Bigg[1+\dfrac{\left(i+\frac{H^2}{\gamma^2}\right)\bigg[i^2p''(\rho)+\kappa_0^2 n^2(1-p'(\rho)) \bigg]\gamma_{w(t)}^2f(w(t))}{2i^2(t)w^2(t)}\Bigg]\dfrac{W}{2ZF(t)}\left(\dfrac{W+Z}{t}\right)\\
    &\hspace{2 cm}+\Bigg[\dfrac{\left(i+\frac{H^2}{\gamma^2}\right)\bigg[i^2p''(\rho)+\kappa_0^2 n^2(1-p'(\rho)) \bigg]\gamma_{w(t)}^2f(w(t))}{2i^2(t)w^2(t)}+2(1-t^2)\Bigg]\dfrac{Wt}{F(t)},\\
    &\hspace{-0.5 cm}Z_t+\dfrac{t^2\sqrt{1-t^2}}{F(t)}Z_r=-\Bigg[1+\dfrac{\left(i+\frac{H^2}{\gamma^2}\right)\bigg[i^2p''(\rho)+\kappa_0^2 n^2(1-p'(\rho)) \bigg]\gamma_{w(t)}^2f(w(t))}{2i^2(t)w^2(t)}\Bigg]\dfrac{Z}{2WF(t)}\left(\dfrac{W+Z}{t}\right)\\
    &\hspace{2 cm}+\Bigg[\dfrac{\left(i+\frac{H^2}{\gamma^2}\right)\bigg[i^2p''(\rho)+\kappa_0^2 n^2(1-p'(\rho)) \bigg]\gamma_{w(t)}^2f(w(t))}{2i^2(t)w^2(t)}+2(1-t^2)\Bigg]\dfrac{Zt}{F(t)}.\\
\end{cases}    
\end{align}
We now derive the boundary information of $W$ and $Z$ in the $(t, r)$ plane using the boundary data \eqref{eq: 3.3} and \eqref{eq: 3.5}. According to the assumption \eqref{eq: 3.6}, it is easy to see that $r=\hat{\theta}(x)$ is strictly decreasing, which means that an inverse function can be defined as $x=\hat{x}(r), r\in [r_1, r_2]$, where $r_1=\hat{\theta}(x_C)$ and $r_2= \hat{\theta}(x_B)$. Therefore the sonic curve $\widehat{BC}: y=\varphi(x)$ can be transformed into 
 a degenerate curve $\widehat{B'C'}: t=0, ~r\in (r_1, r_2]$ in the $(t, r)$ plane such that
\begin{align}
    (W, Z)|_{\widehat{B'C'}}=(-\hat{a}_0, \hat{a}_0)(r), ~r\in (r_1, r_2],
\end{align}
where $\hat{a}_0(r)=\hat{a}_0(\hat{x}(r))$.

Again by differentiating the last equality of \eqref{eq: 3.2} and noticing the assumptions on $\widehat{AB}$, it is obvious that $\tilde{\varpi}'(y)>0$ for all $y\in (y_A, y_B]$, which implies that $\tilde{\theta}'(y)<0$ on $\widehat{BA} / \{B\}$ by \eqref{eq: 3.3a}, which implies the curve $\widehat{BA}: x=\psi(y)$ is transformed into a curve $\widehat{B'A'}: y=\tilde{y}(r),~r\in [r_2, r_3), ~r_3=\tilde{\theta}(y_A)$ such that $W|_{\widehat{B'A'}}=\tilde{b}_0(\tilde{y}(r))$. 

We now proceed to prove that the curve $\widehat{A'B'}$ is actually a $C_+$ characteristic curve in the $(t, r)$ plane:
\begin{l1}
The map $(x, y)\longrightarrow (t,r)$ transforms the positive characteristic curve $\widehat{AB}$ into a positive characteristic curve $\widehat{A'B'}$ of the system \eqref{eq: 4.5} defined by
\begin{align}\label{eq: 4.7}
    r=r_2+\displaystyle \int_0^t \dfrac{t^2\sqrt{1-t^2}}{F(t)} dt=: \bar{r}(t).
\end{align}
\end{l1}
\begin{proof}
Noting that $x(t, r)=\psi(y(t, r))$ on $\widehat{A'B'}$, we differentiate $x(t, r)$ wrt $t$ and use $\psi'=\cot \alpha=\cot(\theta+\omega)$ to yield
  \begin{align*}
      \dfrac{dr}{dt}=\dfrac{y_t\cot (\theta+\omega)-x_t}{x_r-y_r\cot (\theta+\omega)},
  \end{align*}
which along with the fact that $x_t=\dfrac{\theta_y}{J},~x_r=\dfrac{-t_y}{J},~y_t=-\dfrac{\theta_x}{J},~y_r=\dfrac{t_x}{J}$ implies
 \begin{align*}
      \dfrac{dr}{dt}&=-\dfrac{\cot (\theta+\omega)\theta_x+\theta_y}{\sin \omega(\omega_y+\cot (\theta+\omega)\omega_x}\\
      &=-\dfrac{\bar{\partial}_+\theta}{\bar{\partial}_+\omega \sin \omega}=\dfrac{4i^2(t)w^2(t)t^2}{\hat{F}_1(t)\sqrt{1-t^2}}=\dfrac{t^2\sqrt{1-t^2}}{F(t)},
  \end{align*}
which follows that $\widehat{A'B'}$ is a $C_+$ characteristic curve of \eqref{eq: 4.5}. Integrating the above from $0$ to $t$ and noting that $r(0)=r_2$ leads to the expression of \eqref{eq: 4.7} and hence the Lemma is proved. 
\end{proof}
Now noting that 
\begin{align*}
F(0)&=1+\dfrac{\left(i(0)+\frac{H^2(0)}{\gamma^2(0)}\right)\bigg[i^2(0)p''(\rho(0))+\kappa_0^2 n^2(0)(1-p'(\rho(0))\bigg]\gamma_{w(0)}^2f(w(0))}{2i^2(0)w^2(0)}\\
&=\bigg[1+\dfrac{\left(i^2(t)+\frac{H^2}{\gamma^2(t)}\right)\bigg[i^2(t)p''(\rho(t))+\kappa_0^2 n^2(t)(1-p'(\rho(t)) \bigg]\gamma_{w(t)}^2f(w(t))}{2i^2(t)w^2(t)}\bigg]\Bigg|_{t=0}
\end{align*}
and by $W|_{t=0}=-Z|_{t=0},$
we must have 
\begin{align}\label{eq: 4.8}
    W_t|_{t=0}=Z_t|_{t=0}=\dfrac{W+Z}{2t}\bigg|_{t=0}.
\end{align}
In view of the definition of $W$ and $Z$, we see that $\dfrac{W+Z}{2t}=\bar{\partial}_0 I$. Therefore, we must find the boundary value of $\bar{\partial}_0 I$ on $\widehat{BC}$. It follows from \eqref{eq: 2.29} that
\begin{align}\label{eq: 4.9}
    \bar{\partial}_0 I|_{\widehat{BC}}=\dfrac{2i^2(\varpi)w^2(\varpi)}{ F_1(\varpi)} \bar{\partial}_0 \varpi|_{\widehat{BC}},
\end{align}
Again by adding the two equations of \eqref{eq: 2.28}, we have
\begin{align*}
    \bar{\partial}_0 \theta+\dfrac{4i^2(\varpi)w^2(\varpi)\varpi}{F_1(\varpi)} \bar{\partial}_{\perp} \varpi=0.
\end{align*}
Also, $\bar{\partial}_{\perp} \varpi|_{\widehat{BC}}=-\varpi_x(x, \varphi(x))\sin \hat{\theta}+\varpi_y(x, \varphi(x))\cos \hat{\theta}=-F(0)\bar{\partial}_0\theta|_{\widehat{BC}}=-2F(0)\hat{a}_0(x).$
Since $\varpi(x, y)|_{\widehat{BC}}=\hat{\varpi}(x)=1$, we have $\varpi'(x, y)=\varpi_x(x, \varphi(x))+\varphi'(x, \varphi(x))\varpi_y(x, \varphi(x))=0,$ which implies that
\begin{align*}
    \varpi_x(x, \varphi(x))|_{\widehat{BC}}=\dfrac{2\hat{a}_0(x)F(0)\varphi'(x)}{\cos \hat{\theta}+\varphi' \sin \hat{\theta}},~~\varpi_y(x, \varphi(x))|_{\widehat{BC}}=-\dfrac{2\hat{a}_0(x)F(0)}{\cos \hat{\theta}+\varphi' \sin \hat{\theta}}
\end{align*}
or 
\begin{align*}
    \bar{\partial}_0 \varpi|_{\widehat{BC}}=\dfrac{\hat{a}_0(x)F(0)(\varphi' \cos \hat{\theta}-\sin \hat{\theta})}{\cos \hat{\theta}+\varphi' \sin \hat{\theta}}.
\end{align*}
Therefore, \eqref{eq: 4.9} yields
\begin{align}
    \bar{\partial}_0 I|_{\widehat{BC}}=\dfrac{\hat{a}_0(x)(\varphi' \cos \hat{\theta}-\sin \hat{\theta})}{\cos \hat{\theta}+\varphi' \sin \hat{\theta}}=: \hat{a}_1(x) .                             
\end{align}
Let us define $\hat{a}_1(r)=\hat{a}_1(\hat{x}(r))$ and exploit \eqref{eq: 4.8} to obtain the boundary information of $W_t$ and $Z_t$ as $W_t|_{\widehat{B'C'}}=Z_t|_{\widehat{B'C'}}=\hat{a}_1(r)$ such that the function $\hat{a}_1(r)\in C^2((r_1, r_2])$. In summary, the boundary conditions for $(W, Z)$ are as follows:
\begin{align}\label{eq: 4.11}
\begin{cases}
    W(0, r)=-\hat{a}_0(r),~~Z(0, r)=\hat{a}_0(r),~~W_t(0, r)=Z_t(0, r)=\hat{a}_1(r)~~\mathrm{on}~~\mathrm{\widehat{B'C'}},\\
    W(t, \bar{r}(t))=\bar{b}_0(t),~~\mathrm{ on} ~~\mathrm{\widehat{A'B'}},
\end{cases}
\end{align}
where $\bar{b}_0(t)=\tilde{b}_0(\tilde{y}(\bar{r}(t)))$. Also, according to the assumptions \eqref{eq: 3.6}, we know that functions $\hat{a}_0,~\hat{a}_1$ and $\bar{b}_0$ satisfy the following
\begin{align}\label{eq: 4.12}
    \hat{a}_0(r)\in C^2((r_1, r_2]), ~~\hat{a}_1(r)\in C^2((r_1, r_2]),~~\bar{b}_0\in C^2([0, t_0)),~~\hat{a}_0\leq -\epsilon_0, ~~\bar{b}_0\geq \epsilon_0,
\end{align}
where $t_0=\sqrt{1-\tilde{\varpi}^2(A)}$. Moreover, in view of the definitions of $\hat{a}_0, \hat{a}_1$ and $\bar{b}_0$ and the boundary data \eqref{eq: 4.11}, it is easy to see that the compatibility conditions are satisfied at the point $B'$, i.e., $\bar{b}_0(0)=W(0, r_2)=-\hat{a}_0(r_2)$ and $\bar{b}_0'(0)=W_t(0, r_2)=\hat{a}_1(r_2)$. Hence the problem in terms of the $(t, r)$ plane can be reformulated as follows:\\\\
\textbf{Reforumlated problem:} Construct a local classical solution for \eqref{eq: 4.5} with the boundary data \eqref{eq: 4.11} in the region $t>0$ in the neighbourhood of the point $B'(0, r_2)$ under the assumption \eqref{eq: 4.12}.
\subsection{Solution in the partial hodograph plane}
We now construct solution to the reformulated problem in this subsection. Let us denote
\begin{align}\label{eq: 4.13}
    \overline{W}=\dfrac{1}{W}, ~~\overline{Z}=-\dfrac{1}{Z}
\end{align}
to linearize the semi-linear system \eqref{eq: 4.5}. In terms of $\overline{W}$ and $\overline{Z}$, \eqref{eq: 4.5} is transformed into the following linear system
\begin{align}\label{eq: 4.14}
\begin{cases}
    &\hspace{-0.4 cm}\overline{W}_t-\dfrac{t^2\sqrt{1-t^2}}{F(t)}\overline{W}_r=\dfrac{\overline{W}-\overline{Z}}{2t}+tG_1(\overline{W}, \overline{Z}, t),\vspace{0.2 cm}\\
    &\hspace{-0.4 cm}\overline{Z}_t+\dfrac{t^2\sqrt{1-t^2}}{F(t)}\overline{Z}_r=\dfrac{\overline{Z}-\overline{W}}{2t}+tG_2(\overline{W}, \overline{Z}, t),\\
\end{cases}    
\end{align}
where
\begin{align*}
    G_1=&\bigg[\dfrac{\left(i^2(t)+\frac{H^2}{\gamma^2(t)}\right)\bigg[i^2(t)p''(\rho(t))+\kappa_0^2 n^2(t)(1-p'(\rho(t)) \bigg]\gamma_{w(t)}^2f(w(t))}{2i^2(t)w^2(t)}+2-t^2\bigg]\dfrac{(\overline{W}-\overline{Z})}{2F(t)}\\
    &-\bigg[\dfrac{\left(i^2(t)+\frac{H^2}{\gamma^2(t)}\right)\bigg[i^2(t)p''(\rho(t))+\kappa_0^2 n^2(t)(1-p'(\rho(t)) \bigg]\gamma_{w(t)}^2f(w(t))}{2i^2(t)w^2(t)}+2(1-t^2)\bigg]\dfrac{\overline{W}}{F(t)},\\
    G_2=&\bigg[\dfrac{\left(i^2(t)+\frac{H^2}{\gamma^2(t)}\right)\bigg[i^2(t)p''(\rho(t))+\kappa_0^2 n^2(t)(1-p'(\rho(t)) \bigg]\gamma_{w(t)}^2f(w(t))}{2i^2(t)w^2(t)}+2-t^2\bigg]\dfrac{(\overline{Z}-\overline{W})}{2F(t)}\\
    &-\bigg[\dfrac{\left(i^2(t)+\frac{H^2}{\gamma^2(t)}\right)\bigg[i^2(t)p''(\rho(t))+\kappa_0^2 n^2(t)(1-p'(\rho(t)) \bigg]\gamma_{w(t)}^2f(w(t))}{2i^2(t)w^2(t)}+2(1-t^2)\bigg]\dfrac{\overline{Z}}{F(t)}.\\
\end{align*}
We now use the transformation $\upsilon=t, ~~\chi=\bar{r}(t)-r$ to reduce \eqref{eq: 4.14} into the following form
\begin{align}\label{eq: 4.15}
    \begin{cases}
    &\hspace{-0.4 cm}\widetilde{W}_{\upsilon}+\dfrac{2\upsilon^2\sqrt{1-\upsilon^2}}{F(\upsilon)}\widetilde{W}_\chi=\dfrac{\widetilde{W}-\widetilde{Z}}{2\upsilon}+\upsilon G_1(\widetilde{W}, \widetilde{Z}, \upsilon),\vspace{0.2 cm}\\
    &\hspace{-0.4 cm}\widetilde{Z}_{\upsilon}=\dfrac{\widetilde{Z}-\widetilde{W}}{2\upsilon}+\upsilon G_2(\widetilde{W}, \widetilde{Z}, \upsilon),
\end{cases} 
\end{align}
where $\widetilde{W}(\upsilon, \chi)=\overline{W}(\upsilon, \bar{r}(t)-\chi)$ and $\widetilde{Z}(\upsilon, \chi)=\overline{Z}(\upsilon, \bar{r}(t)-\chi)$. Recalling the boundary data \eqref{eq: 4.11}, we now have the boundary data of $\widetilde{W}(\upsilon, \chi), \widetilde{Z}(\upsilon, \chi)$ as follows:
\begin{align}\label{eq: 4.16}
\begin{cases}
    \widetilde{W}(0, \chi)=\widetilde{Z}(0, \chi)=a_0(\chi), ~~\widetilde{W}_{\upsilon}(0, \chi)=-\widetilde{Z}_{\upsilon}(0, \chi)=a_1(\chi)~~~~\mathrm{on}~~\upsilon=0, ~~0\leq \chi<r_2-r_1,\\
    \widetilde{W}(\chi, 0)=b_0(\upsilon)~~~~\mathrm{on}~~\chi=0, ~~0\leq \upsilon<t_0,
\end{cases}    
\end{align}
where 
\begin{align}
    a_0(\chi)=-\dfrac{1}{\hat{a}_0(r_2-\chi)}, ~~ a_1(\chi)=\dfrac{\hat{a}_1}{\hat{a}_0^2}(r_2-\chi),~~ b_0(\upsilon)=\dfrac{1}{\bar{b}_0(\upsilon)}.
\end{align}
We now denote
\begin{align}
    U=\widetilde{W}-a_0(\chi)+a_1(\chi)\upsilon, ~~~~~V=\widetilde{Z}-a_0(\chi)-a_1(\chi)\upsilon
\end{align}
to homogenize the boundary conditions \eqref{eq: 4.11} such that
\begin{align}\label{eq: 4.19}
\begin{cases}
    U(0, \chi)=V(0, \chi)=U_{\upsilon}(0, \chi)=V_{\upsilon}(0, \chi)=0, \forall \chi~\in [0, r_2-r_1),\\
    U(\upsilon, 0)=b_1(\upsilon), ~\forall~ \upsilon\in [0, t_0),
\end{cases}    
\end{align}
where the function $b_1(\upsilon)$ satisfies $b_1(\upsilon)=b_0(\upsilon)-a_0(0)+a_1(0)\upsilon$, which implies that $b_1'(\upsilon)=b_0'(\upsilon)+a_1(0)$. Therefore, $b_1'(0)=b_1(0)=0$.

By a direct calculation, one can convert the system \eqref{eq: 4.15} in terms of $(U, V)$ as follows:
\begin{align}\label{eq: 4.20}
    \begin{cases}
    &\hspace{-0.4 cm}U_{\upsilon}+\dfrac{2\upsilon^2\sqrt{1-\upsilon^2}}{F(\upsilon)}U_\chi=\dfrac{U-V}{2\upsilon}+\Bigg[\dfrac{\left(i+\frac{H^2}{\gamma^2}\right)\bigg[ i^2p''(\rho)+\kappa_0^2 n^2(1-p'(\rho))\bigg]\gamma_{w}^2f(w)}{2i^2w^2}+2-\upsilon^2\Bigg]\dfrac{(U-V)\upsilon }{2F(\upsilon)}\\
    &\hspace{1 cm}-\Bigg[\dfrac{\left(i+\frac{H^2}{\gamma^2}\right)\bigg[ i^2p''(\rho)+\kappa_0^2 n^2(1-p'(\rho))\bigg]\gamma_{w}^2f(w)}{2i^2w^2}+2(1-\upsilon^2)\Bigg]\dfrac{U\upsilon}{F(\upsilon)}+E_1(\upsilon, \chi)\upsilon,\vspace{0.2 cm}\\
    &\hspace{-0.4 cm}V_{\upsilon}=\dfrac{V-U}{2\upsilon}+\Bigg[\dfrac{\left(i+\frac{H^2}{\gamma^2}\right)\bigg[ i^2p''(\rho)+\kappa_0^2 n^2(1-p'(\rho))\bigg]\gamma_{w}^2f(w)}{2i^2w^2}+2-\upsilon^2\Bigg]\dfrac{(V-U)\upsilon }{2F(\upsilon)}\\
    &\hspace{1 cm}-\Bigg[\dfrac{\left(i+\frac{H^2}{\gamma^2}\right)\bigg[ i^2p''(\rho)+\kappa_0^2 n^2(1-p'(\rho))\bigg]\gamma_{w}^2f(w)}{2i^2w^2}+2(1-\upsilon^2)\Bigg]\dfrac{V\upsilon}{F(\upsilon)}+E_2(\upsilon, \chi)\upsilon,
\end{cases} 
\end{align}
where 
\begin{align*}
    E_1(\upsilon, \chi)&=-\bigg[\dfrac{\left(i+\frac{H^2}{\gamma^2}\right)\left(i^2p''(\rho)+\kappa_0^2 n^2(1-p'(\rho)) \right)\gamma_{w(\upsilon)}^2f(w(\upsilon))}{2i^2(\upsilon)w^2(\upsilon)}+2(1-\upsilon^2)\bigg]\dfrac{a_0}{F(\upsilon)}\\
    &-\dfrac{2\upsilon \sqrt{1-\upsilon^2}}{F}(a_0'-a_1' \upsilon)-\dfrac{a_1\upsilon^3}{F(\upsilon)}
\end{align*}
and 
\begin{align*}
    E_2(\upsilon, \chi)&=-\Bigg[\dfrac{\left(i+\frac{H^2}{\gamma^2}\right)\bigg[ i^2p''(\rho)+\kappa_0^2 n^2(1-p'(\rho))\bigg]\gamma_{w(\upsilon)}^2f(w(\upsilon))}{2i^2(\upsilon)w^2(\upsilon)}+2(1-\upsilon^2)\Bigg]a_0+\dfrac{a_1\upsilon^3}{F(\upsilon)}.
\end{align*}
In view of the regularity of $a_0(\chi), a_1(\chi)$, it is easy to see that $F_1$ and $F_2$ have continuous derivatives wrt $\chi$. Therefore the problem is converted to seek a classical solution for a degenerate linear initial-characteristic problem, which has two eigenvalues given by
\begin{align}
    \lambda_-=0, ~~\lambda_+=\dfrac{2\upsilon^2\sqrt{1-\upsilon^2}}{F(\upsilon)}.
\end{align}
We now denote $\chi=\bar{\chi}(\upsilon)$ as the $C_+$ characteristic passing through the origin, which can be given by
\begin{align*}
    \bar{\chi}(\upsilon)=\displaystyle \int_0^{\upsilon} \dfrac{2s^2\sqrt{1-s^2}}{F(s)} ds.
\end{align*}
Further, we define a square region $\Omega:= \{(\upsilon, \chi)| 0\leq \upsilon \leq \delta, ~0\leq \chi\leq \delta\}$ for a positive constant $\delta$. Then we integrate the system \eqref{eq: 4.20} with the boundary conditions \eqref{eq: 4.19} along the characteristics to yield a system of integral equations given by
{
\footnotesize
\begin{align}\label{eq: 4.22}
\begin{cases}
V(\tau, \zeta)&=\displaystyle \bigintssss_0^{\tau} \Bigg\{ \dfrac{V-U}{2\upsilon}+\Bigg[\dfrac{\left(i+\frac{H^2}{\gamma^2}\right)\bigg[ i^2p''(\rho)+\kappa_0^2 n^2(1-p'(\rho))\bigg]\gamma_{w(\upsilon)}^2f(w(\upsilon))}{2i^2(\upsilon)w^2(\upsilon)}+2-\upsilon^2\Bigg]\dfrac{(V-U)\upsilon }{2F}\\
&\hspace{2 cm}-\Bigg[\dfrac{\left(i+\frac{H^2}{\gamma^2}\right)\bigg[ i^2p''(\rho)+\kappa_0^2 n^2(1-p'(\rho))\bigg]\gamma_{w(\upsilon)}^2f(w(\upsilon))}{2i^2(\upsilon)w^2(\upsilon)}+2(1-\upsilon^2)\Bigg]\dfrac{V\upsilon}{F}+E_2\upsilon\Bigg\}(\upsilon, \chi_-(\upsilon)) d\upsilon,\vspace{0.3 cm}\\
U(\tau, \zeta)&=\begin{cases}
    \displaystyle \bigintssss_0^{\tau} \Bigg\{\dfrac{U-V}{2\upsilon}+\Bigg[\dfrac{\left(i+\frac{H^2}{\gamma^2}\right)\bigg[ i^2p''(\rho)+\kappa_0^2 n^2(1-p'(\rho))\bigg]\gamma_{w(\upsilon)}^2f(w(\upsilon))}{2i^2(\upsilon)w^2(\upsilon)}+2-\upsilon^2\Bigg]\dfrac{(U-V)\upsilon }{2F}\\
    \hspace{0.5 cm}-\Bigg[\dfrac{\left(i+\frac{H^2}{\gamma^2}\right)\bigg[ i^2p''(\rho)+\kappa_0^2 n^2(1-p'(\rho))\bigg]\gamma_{w(\upsilon)}^2f(w(\upsilon))}{2i^2(\upsilon)w^2(\upsilon)}+2(1-\upsilon^2)\Bigg]\dfrac{U\upsilon}{F(\upsilon)}+E_1\upsilon\Bigg\}(\upsilon, \chi_+(\upsilon)) d\upsilon,\\
    \hspace{12 cm}\zeta\geq \bar{\chi}(\tau), \tau \geq 0,\vspace{0.3 cm}\\
    b_1(\tau_1)+\displaystyle \bigintssss_{\tau_1}^{\tau} \Bigg\{\dfrac{U-V}{2\upsilon}+\Bigg[\dfrac{\left(i+\frac{H^2}{\gamma^2}\right)\bigg[ i^2p''(\rho)+\kappa_0^2 n^2(1-p'(\rho))\bigg]\gamma_{w(\upsilon)}^2f(w(\upsilon))}{2i^2(\upsilon)w^2(\upsilon)}+2-\upsilon^2\Bigg]\dfrac{(U-V)\upsilon }{2F}\\
    \hspace{0.5 cm}-\Bigg[\dfrac{\left(i+\frac{H^2}{\gamma^2}\right)\bigg[ i^2p''(\rho)+\kappa_0^2 n^2(1-p'(\rho))\bigg]\gamma_{w(\upsilon)}^2f(w(\upsilon))}{2i^2(\upsilon)w^2(\upsilon)}+2(1-\upsilon^2)\Bigg]\dfrac{U\upsilon}{F(\upsilon)}+E_1\upsilon\Bigg\}(\upsilon, \chi_+(\upsilon)) d\upsilon\\
    \hspace{12 cm}0\leq \zeta<\bar{\chi}(\tau), \tau \geq 0,
\end{cases} 
\end{cases}
\end{align}
}
where $\chi_+(\upsilon)=\chi_+(\upsilon; \tau, \zeta), \chi_-(\upsilon)=\chi_-(\upsilon; \tau, \zeta)$ are the positive and negative characteristics passing through a point $(\tau, \zeta)$, respectively and $\tau_1\in [0, \tau)$ is a function of $(\tau, \zeta)$ defined as
\begin{align}\label{eq: 4.23}
    \zeta=\displaystyle \int_{\tau_1}^{\tau} \dfrac{2\upsilon^2\sqrt{1-\upsilon^2}}{F(\upsilon)} d\upsilon.
\end{align}
Moreover, by using the fact that $\lambda_-=0$ we must have $\chi_-(\upsilon; \tau, \zeta)=\zeta$.

We now prove the existence of solution for the system \eqref{eq: 4.22} by exploiting the method of iterations. We denote $U^0(\upsilon, \chi)=V^0(\upsilon, \chi)=0$ and define $U^{(k)}$ and $V^{(k)}, (k\geq 1)$ as follows:
{
\footnotesize
\begin{align}\label{eq: 4.24}
\begin{cases}
\hspace{-0.1 cm}V^{(k)}(\tau, \zeta)&=\displaystyle \bigintssss_0^{\tau} \Bigg\{ \dfrac{V^{(k-1)}-U^{(k-1)}}{2\upsilon} -\Bigg[\dfrac{\left(i+\frac{H^2}{\gamma^2}\right)\bigg[ i^2p''(\rho)+\kappa_0^2 n^2(1-p'(\rho))\bigg]\gamma_{w(\upsilon)}^2f(w(\upsilon))}{2i^2(\upsilon)w^2(\upsilon)}+2(1-\upsilon^2)\Bigg]\dfrac{V^{(k-1)}\upsilon}{F}\\
&+\Bigg[\dfrac{\left(i+\frac{H^2}{\gamma^2}\right)\bigg[ i^2p''(\rho)+\kappa_0^2 n^2(1-p'(\rho))\bigg]\gamma_{w(\upsilon)}^2f(w(\upsilon))}{2i^2(\upsilon)w^2(\upsilon)}+2-\upsilon^2\Bigg]\dfrac{(V^{(k-1)}-U^{(k-1)})\upsilon}{2F}+E_2\upsilon\Bigg\}(\upsilon, \chi_-(\upsilon)) d\upsilon,\vspace{0.3 cm}\\
\hspace{-0.1 cm}U^{(k)}(\tau, \zeta)\hspace{-0.1 cm}&=\begin{cases}
    \hspace{-0.1 cm}\displaystyle \bigintssss_0^{\tau} \Bigg\{\dfrac{U^{(k-1)}-V^{(k-1)}}{2\upsilon}-\Bigg[\dfrac{\left(i+\frac{H^2}{\gamma^2}\right)\bigg[ i^2p''(\rho)+\kappa_0^2 n^2(1-p'(\rho))\bigg]\gamma_{w(\upsilon)}^2f(w(\upsilon))}{2i^2(\upsilon)w^2(\upsilon)}+2(1-\upsilon^2)\Bigg]\dfrac{U^{(k-1)}\upsilon}{F(\upsilon)}\\
  \hspace{-0.1 cm}+\Bigg[\dfrac{\left(i+\frac{H^2}{\gamma^2}\right)\bigg[ i^2p''(\rho)+\kappa_0^2 n^2(1-p'(\rho))\bigg]\gamma_{w}^2f(w)}{2i^2(\upsilon)w^2(\upsilon)}+2-\upsilon^2\Bigg]\dfrac{(U^{(k-1)}-V^{(k-1)})\upsilon }{2F}+E_1\upsilon\Bigg\}(\upsilon, \chi_+(\upsilon)) d\upsilon,\\
    \hspace{12 cm}\zeta\geq \bar{\chi}(\tau), \tau \geq 0,\vspace{0.3 cm}\\
  \hspace{-0.1 cm}\displaystyle \bigintssss_{\tau_1}^{\tau} \Bigg\{\dfrac{U^{(k-1)}-V^{(k-1)}}{2\upsilon}-\Bigg[\dfrac{\left(i+\frac{H^2}{\gamma^2}\right)\bigg[ i^2p''(\rho)+\kappa_0^2 n^2(1-p'(\rho))\bigg]\gamma_{w}^2f(w)}{2i^2(\upsilon)w^2(\upsilon)}+2(1-\upsilon^2)\Bigg]\dfrac{U^{(k-1)}\upsilon}{F(\upsilon)}\\
  \hspace{-0.1 cm}+\Bigg[\dfrac{\left(i+\frac{H^2}{\gamma^2}\right)\bigg[ i^2p''(\rho)+\kappa_0^2 n^2(1-p'(\rho))\bigg]\gamma_{w}^2f(w)}{2i^2(\upsilon)w^2(\upsilon)}+2-\upsilon^2\Bigg]\dfrac{(U^{(k-1)}-V^{(k-1)})\upsilon }{2F}+E_1\upsilon\Bigg\}(\upsilon, \chi_+(\upsilon)) d\upsilon\\
    +  b_1(\tau_1),\hspace{11 cm}0\leq \zeta<\bar{\chi}(\tau), \tau \geq 0.
\end{cases} 
\end{cases}
\end{align}
}
Now we proceed to prove that the sequences $\{(U^{(k)}, V^{(k)})\}$ are uniformly convergent in the domain $\Omega$ for a small positive constant $\delta.$ In view of the expression of $F(t)$ and noting that $0<a^2<1$, $0<q<\hat{q}<1$, it is straightforward to observe that there exists a sufficiently small $\delta>0$ and $\upsilon\in [0, \delta]$ such that $0<\dfrac{\sqrt{1-\upsilon^2}}{F(\upsilon)}\leq \dfrac{1}{F(\upsilon)}\leq K_\delta<\infty$, where $K_\delta$ is a uniform small positive constant. Moreover, we must have
\begin{align}
   0<\dfrac{1}{F(\upsilon)}\left(\dfrac{\left(i+\frac{H^2}{\gamma^2}\right)\bigg[ i^2p''(\rho)+\kappa_0^2 n^2(1-p'(\rho))\bigg]\gamma_{w}^2f(w)}{2i^2w^2}+2(1-\upsilon^2)\right)\leq K_\delta,\\
    0<\dfrac{1}{F(\upsilon)}\left(\dfrac{\left(i+\frac{H^2}{\gamma^2}\right)\bigg[ i^2p''(\rho)+\kappa_0^2 n^2(1-p'(\rho))\bigg]\gamma_{w}^2f(w)}{2i^2w^2}+2-\upsilon^2\right)\leq K_\delta,
\end{align}
For simplicity of notations, let $M>0$ be a uniform constant depending only on the $C^2$ norms of $a_0, a_1$ and $b_0$ which can change between different expressions in the analysis of this article. From the expressions of $F_1$ and $F_2$, it is clear that
\begin{align}
    |E_1|+|E_{1\chi}|+|E_2|+|E_{2\chi}|\leq M.
\end{align}
For any point $P(\tau, \zeta)$ in $\Omega=\{(\upsilon, \chi)| 0< \upsilon \leq \delta, ~0<\chi\leq \delta\}$, by the expression of $b_1(.)$ it is clear that $|b_1(\tau_1)|=|b_0(\tau_1)-a_0(0)+a_1(0)\tau_1|\leq \dfrac{M}{2}\tau_1^2$ if $\zeta<\bar{\chi}(\tau)$. Noting that $\tau_1\leq \tau$, we must have
\begin{align}\label{eq: 4.28}
    |b_1(\tau_1)|\leq \dfrac{M\tau^2}{2}.
\end{align}
Further, employing \eqref{eq: 4.23} for any $\bar{\tau}\in [0, \tau]$, we must have
\begin{align}
    |\chi_+(\bar{\tau})-\chi_-(\bar{\tau})|&= |\chi_+(\bar{\tau}; \tau, \zeta)-\chi_-(\bar{\tau}; \tau, \zeta)| \nonumber\\
    & \leq |\chi_+(0; \tau, \zeta)- \zeta|=\bigg|\displaystyle \int_{0}^{\tau} \dfrac{2\upsilon^2 \sqrt{1-\upsilon^2}}{F(\upsilon)}d\upsilon\bigg|\leq \dfrac{2K_\delta \tau^3}{3}~~\mathrm{if}~~ \zeta\geq \bar{\chi}(\tau)
\end{align}
and for any $\bar{\tau}\in [\bar{\tau}_1, \tau]$, we have
\begin{align}
    |\chi_+(\bar{\tau})-\chi_-(\bar{\tau})|&= |\chi_+(\bar{\tau}; \tau, \zeta)-\chi_-(\bar{\tau}; \tau, \zeta)|\nonumber \\
    & \leq |\chi_+(\bar{\tau_1}; \tau, \zeta)- \zeta|=\bigg|\displaystyle \int_{\tau_1}^{\tau} \dfrac{2\upsilon^2 \sqrt{1-\upsilon^2}}{F(\upsilon)}d\upsilon\bigg|\leq \dfrac{2K_\delta \tau^3}{3}~~\mathrm{if}~~ \zeta< \bar{\chi}(\tau).
\end{align}
Based on the above discussions, one can derive the following Lemmas:
\begin{l1}\label{l-4.2}
For any $k\geq 1$ and $(\tau, \zeta)\in \Omega$, the following inequalities hold true
\begin{align}
\begin{cases}
    |U^{(k)}(\tau, \zeta)|; |V^{(k)}(\tau, \zeta)|\leq M\tau^2 \sum_{j=0}^{k} \left(\dfrac{2}{3}\right)^j,\\
    |U^{(k)}(\tau, \zeta)-V^{(k)}(\tau, \zeta)|\leq M\tau^2 \sum_{j=0}^{k} \left(\dfrac{2}{3}\right)^j.
\end{cases}
\end{align}
\end{l1}
\begin{proof}
We use the induction argument to prove this Lemma. The proof can be discussed into two different cases depending on whether $\zeta<\bar{\chi}(\tau)$ or $\zeta \geq \bar{\chi}(\tau)$. In this proof, we analyze the first case $\zeta<\bar{\chi}(\tau)$ only as the second case can be analyzed similarly.

Now for $n=1$, we have 
\begin{align}
    |V^{(1)}(\tau, \zeta)|\leq \displaystyle \int_0^{\tau}|\upsilon E_2(\upsilon)| d\upsilon\leq \displaystyle \int_0^{\tau}M\upsilon d\upsilon \leq  \dfrac{M\tau^2}{2}\leq M\tau^2 \sum_{j=0}^{1} \left(\dfrac{2}{3}\right)^j.
\end{align}
Similarly we estimate $U^{(1)}$ for $\zeta<\bar{\chi}(\tau)$ as
\begin{align}
    |U^{(1)}(\tau, \zeta)|&\leq \displaystyle \int_{\tau_1}^{\tau}|\upsilon E_1(\upsilon)| d\upsilon+|b_1(\tau_1)|\leq \displaystyle \int_{\tau_1}^{\tau}M\upsilon d\upsilon+\dfrac{M\tau^2}{2} \nonumber \\
    &\leq \dfrac{M(\tau^2-\tau_1^2)}{2}+\dfrac{M\tau^2}{2}\leq M\tau^2 \sum_{j=0}^{1} \left(\dfrac{2}{3}\right)^j.
\end{align}
Now we proceed to estimate $|U^{(1)}(\tau, \zeta)-V^{(1)}(\tau, \zeta)|$. Now by direct calculation, one can compute that
\begin{align}\label{eq: 4.35}
|U^{(1)}(\tau, \zeta)&-V^{(1)}(\tau, \zeta)| \nonumber\\
&\leq \displaystyle \int_{\tau_1}^{\tau}|E_1(\upsilon, \chi_+(\upsilon; \tau, \zeta))-E_2(\upsilon, \chi_-(\upsilon; \tau, \zeta))| \upsilon d\upsilon+\displaystyle \int_0^{\tau_1}|E_2(\upsilon, \chi_-(\upsilon; \tau, \zeta))| \upsilon d\upsilon+|b_1(\tau_1)|.
\end{align}
Hence, it is clear that we need to estimate $|E_1(\upsilon, \chi_+(\upsilon; \tau, \zeta))-E_2(\upsilon, \chi_-(\upsilon; \tau, \zeta))|$ in order to estimate $|U^{(1)}(\tau, \zeta)-V^{(1)}(\tau, \zeta)|$. Now
\begin{align*}
|E_1(\upsilon,& \chi_+(\upsilon; \tau, \zeta))-E_2(\upsilon, \chi_-(\upsilon; \tau, \zeta))| \nonumber\\
&\leq \dfrac{2\upsilon \sqrt{1-\upsilon^2}}{F}(|a_0'|+|a_1'|\upsilon)+\dfrac{|a_1(\chi_+(\upsilon; \tau, \zeta))+a_1(\chi_-(\upsilon; \tau, \zeta))|}{F}\upsilon^3\\
&\hspace{-0.5 cm}+\Bigg[\dfrac{\left(i+\frac{H^2}{\gamma^2}\right)\bigg[ i^2p''(\rho)+\kappa_0^2 n^2(1-p'(\rho))\bigg]\gamma_{w}^2f(w)}{2i^2(\upsilon)w^2(\upsilon)}+2(1-\upsilon^2)\Bigg]\dfrac{|a_0(\chi_+(\upsilon; \tau, \zeta))-a_0(\chi_-(\upsilon; \tau, \zeta))|}{F} \nonumber\\
& \leq 2 K_\delta(M+M\upsilon)\upsilon +2MK_\delta \upsilon^3+K_\delta M|\chi_+-\chi_-| \nonumber\\
& \leq 2K_\delta (M+M\tau)\tau+2MK_\delta \tau^3+\dfrac{2}{3}MK_\delta \tau^3\\
& \leq K_\delta M\tau,
\end{align*}
which is exploited in \eqref{eq: 4.35} along with \eqref{eq: 4.28} to obtain
\begin{align}
|U^{(1)}(\tau, \zeta)-V^{(1)}(\tau, \zeta)|& \leq \displaystyle \int_{\tau_1}^{\tau} MK_\delta \tau \upsilon d\upsilon+\displaystyle \int_0^{\tau_1} M\upsilon d\upsilon+\dfrac{1}{2}M\tau^2 \nonumber\\
&\leq \dfrac{MK_{\delta}}{2}\tau(\tau^2-\tau_1^2)+\dfrac{1}{2}M\tau_1^2+\dfrac{1}{2}M\tau^2 \nonumber\\
&\leq \dfrac{MK_{\delta}}{2}\tau^3+\dfrac{1}{2}M\tau^2+\dfrac{1}{2}M\tau^2 \leq M\tau^2 \left(\dfrac{\delta K_{\delta}}{2}+1\right) \nonumber\\
&\leq M\tau^2 \sum_{j=0}^{1} \left(\dfrac{2}{3}\right)^j.
\end{align}
Now we assume that the induction argument is true for $n=k$, and then we need to prove that the statement is true for $n=k+1$. Now we use the induction assumption to compute
\begin{align*}
    |V^{(k+1)}(\tau, \zeta)|&\leq \displaystyle \bigintssss_{0}^{\tau} \Bigg|\dfrac{V^{(k)}(\upsilon, \chi_-(\upsilon))-U^{(k)}(\upsilon, \chi_-(\upsilon))}{2\upsilon}+E_2(\upsilon, \chi_-(\upsilon))\upsilon\Bigg| d\upsilon\\
    &\hspace{-1 cm}+\displaystyle \bigintssss_{0}^{\tau} \Bigg|\Bigg[\dfrac{\left(i+\frac{H^2}{\gamma^2}\right)\bigg[ i^2p''(\rho)+\kappa_0^2 n^2(1-p'(\rho))\bigg]\gamma_{w(\upsilon)}^2f(w(\upsilon))}{2i^2(\upsilon)w^2(\upsilon)}+2(1-\upsilon^2)\Bigg]\dfrac{V^{(k)}(\upsilon, \chi_-(\upsilon))\upsilon}{F}\Bigg|d\upsilon\\
    &\hspace{-3 cm}+\displaystyle \bigintssss_{0}^{\tau} \Bigg|\Bigg[\dfrac{\left(i+\frac{H^2}{\gamma^2}\right)\bigg[ i^2p''(\rho)+\kappa_0^2 n^2(1-p'(\rho))\bigg]\gamma_{w}^2f(w)}{2i^2(\upsilon)w^2(\upsilon)}+2-\upsilon^2\Bigg]\dfrac{\left(U^{(k)}(\upsilon, \chi_-(\upsilon))-V^{(k)}(\upsilon, \chi_-(\upsilon))\right)\upsilon}{2F}\bigg| d\upsilon\\
    & \leq \displaystyle \int_{0}^{\tau} \bigg\{\bigg[\dfrac{M\upsilon}{2}+2K_\delta M\upsilon^3\bigg]\sum_{j=0}^{k} \left(\dfrac{2}{3}\right)^j+M\upsilon\bigg\}d\upsilon= M\tau^2\bigg\{\dfrac{1}{2}+\bigg[\dfrac{1}{4}+\dfrac{\delta^2K_\delta}{2}\bigg]\sum_{j=0}^{k} \left(\dfrac{2}{3}\right)^j\bigg\}\\
    &\leq M\tau^2 \sum_{j=0}^{k+1} \left(\dfrac{2}{3}\right)^j.
\end{align*}
One can find the following estimate of $U^{(k+1)}$ in a similar manner
\begin{align}
    |U^{(k+1)}(\tau, \zeta)|\leq  M\tau^2 \sum_{j=0}^{k+1} \left(\dfrac{2}{3}\right)^j.
\end{align}
Now we compute the estimates of $|U^{(k+1)}(\tau, \zeta)-V^{(k+1)}(\tau, \zeta)|$. A straightforward calculation yields
\begin{align}
    |U^{(k+1)}(\tau, \zeta)-V^{(k+1)}(\tau, \zeta)|\leq I_1+I_2+|b_1(\tau_1)|,
\end{align}
where 
\begin{align*}
    I_1&=\displaystyle \bigintssss_{\tau_1}^{\tau} \Bigg\{ \dfrac{|(U^{(k)}-V^{(k)})|(\upsilon, \chi_+(\upsilon; \tau, \zeta))-|U^{(k)}-V^{(k)}|(\upsilon, \chi_-(\upsilon; \tau, \zeta))}{2\upsilon}\\
    &+\left(\dfrac{\left(i+\frac{H^2}{\gamma^2}\right)\bigg[ i^2p''(\rho)+\kappa_0^2 n^2(1-p'(\rho))\bigg]\gamma_{w(\upsilon)}^2f(w(\upsilon))}{2i^2(\upsilon)w^2(\upsilon)}+2-\upsilon^2\right)\\
    &\hspace{4 cm}\dfrac{\bigg[|U^{(k)}-V^{(k)}|(\upsilon, \chi_+(\upsilon; \tau, \zeta))+|U^{(k)}-V^{(k)}|(\upsilon, \chi_-(\upsilon; \tau, \zeta))\bigg]\upsilon}{2F}\\
    &\hspace{-2 cm}+\left(\dfrac{\left(i+\frac{H^2}{\gamma^2}\right)\bigg[ i^2p''(\rho)+\kappa_0^2 n^2(1-p'(\rho))\bigg]\gamma_{w}^2f(w)}{2i^2(\upsilon)w^2(\upsilon)}+2(1-\upsilon^2)\right)\dfrac{\bigg[|U^{(k)}(\upsilon, \chi_+(\upsilon; \tau, \zeta))|+|V^{(k)}(\upsilon, \chi_-(\upsilon; \tau, \zeta))|\bigg]\upsilon}{F}\\
    &\hspace{1 cm}+|E_1(\upsilon, \chi_+(\upsilon; \tau, \zeta))-E_2(\upsilon, \chi_-(\upsilon; \tau, \zeta))|\upsilon\Bigg\} d\upsilon
\end{align*}
and
\begin{align*}
    I_2&=\displaystyle \bigintssss_{0}^{\tau_1} \Bigg\{ \dfrac{|(U^{(k)}-V^{(k)})(\upsilon, \chi_-(\upsilon; \tau, \zeta))|}{2\upsilon}+|E_2(\upsilon, \chi_-(\upsilon; \tau, \zeta))|\upsilon\\
    &+\left(\dfrac{\left(i+\frac{H^2}{\gamma^2}\right)\bigg[ i^2p''(\rho)+\kappa_0^2 n^2(1-p'(\rho))\bigg]\gamma_{w(\upsilon)}^2f(w(\upsilon))}{2i^2(\upsilon)w^2(\upsilon)}+2-\upsilon^2\right)\dfrac{|(U^{(k)}-V^{(k)})(\upsilon, \chi_-(\upsilon; \tau, \zeta))|\upsilon}{2F}\\
    &+\left(\dfrac{\left(i+\frac{H^2}{\gamma^2}\right)\bigg[ i^2p''(\rho)+\kappa_0^2 n^2(1-p'(\rho))\bigg]\gamma_{w(\upsilon)}^2f(w(\upsilon))}{2i^2(\upsilon)w^2(\upsilon)}+2(1-\upsilon^2)\right)\dfrac{|V^{(k)}(\upsilon, \chi_-(\upsilon; \tau, \zeta))|\upsilon}{F}\Bigg\} d\upsilon.
\end{align*}
Clearly, we need to find the estimates of $I_1$ and $I_2$ in order to find the estimate of $|U^{(k+1)}(\tau, \zeta)-V^{(k+1)}(\tau, \zeta)|$. 

Now using the induction assumption, we have 
\begin{align*}
    I_1&\leq \displaystyle \int_{\tau_1}^{\tau} \bigg\{M\upsilon \sum_{j=0}^{k} \left(\dfrac{2}{3}\right)^j+2K_\delta M\upsilon^3 \sum_{j=0}^{k} \left(\dfrac{2}{3}\right)^j+K_\delta M \upsilon^3 \sum_{j=0}^{k} \left(\dfrac{2}{3}\right)^j+K_\delta M \tau \upsilon \bigg\}d\upsilon\\
    & \leq \dfrac{M}{2}(\tau^2-\tau_1^2)\sum_{j=0}^{k} \left(\dfrac{2}{3}\right)^j+\dfrac{3K_\delta M}{4}(\tau^4-\tau_1^4)\sum_{j=0}^{k} \left(\dfrac{2}{3}\right)^j+\dfrac{K_\delta M\tau^3}{2},
\end{align*}
and 
\begin{align*}
    I_2&\leq \displaystyle \int_{0}^{\tau_1} \bigg\{\dfrac{M\upsilon}{2} \sum_{j=0}^{k} \left(\dfrac{2}{3}\right)^j+2K_\delta M\upsilon^3 \sum_{j=0}^{k} \left(\dfrac{2}{3}\right)^j+ M \upsilon \bigg\}d\upsilon\\
    & \leq \dfrac{M\tau_1^2}{4}\sum_{j=0}^{k} \left(\dfrac{2}{3}\right)^j+\dfrac{K_\delta M\tau_1^4}{2}\sum_{j=0}^{k} \left(\dfrac{2}{3}\right)^j+\dfrac{M\tau_1^2}{2},
\end{align*}
which implies 
\begin{align*}
    &|U^{(k+1)}(\tau, \zeta)-V^{(k+1)}(\tau, \zeta)|\\
    &\leq \bigg\{\dfrac{M}{2}(\tau^2-\tau_1^2)\sum_{j=0}^{k} \left(\dfrac{2}{3}\right)^j+\dfrac{3K_\delta M}{4}(\tau^4-\tau_1^4)\sum_{j=0}^{k} \left(\dfrac{2}{3}\right)^j+\dfrac{K_\delta M\tau^3}{2}\bigg\}\\
    &~~~+\bigg\{\dfrac{M\tau_1^2}{4}\sum_{j=0}^{k} \left(\dfrac{2}{3}\right)^j+\dfrac{K_\delta M\tau_1^4}{2}\sum_{j=0}^{k} \left(\dfrac{2}{3}\right)^j+\dfrac{M\tau_1^2}{2}\bigg\}+\dfrac{M}{2}\tau^2\\
    & \leq M\tau^2 \bigg\{ \dfrac{1}{2}+\dfrac{\delta K_\delta}{2}+\bigg[\dfrac{3\delta^2 K_\delta}{4}+\dfrac{1}{2}\bigg]\sum_{j=0}^{k} \left(\dfrac{2}{3}\right)^j\bigg\}\leq M\tau^2  \sum_{j=0}^{k+1} \left(\dfrac{2}{3}\right)^j,
\end{align*}
which shows that the proof of the Lemma is complete.
\end{proof}
\begin{l1}\label{l-4.3}
For all $k\geq 1$ and any $(\tau, \zeta)\in \Omega$, the following inequalities hold true
\begin{align}\label{eq: 4.38}
    \begin{cases}
    |U^{(k+1)}(\tau, \zeta)-U^{(k)}(\tau, \zeta)|; |V^{(k+1)}(\tau, \zeta)-V^{(k)}(\tau, \zeta)|\leq M\tau^2  \left(\dfrac{2}{3}\right)^j,\\
    |U^{(k+1)}(\tau, \zeta)-V^{(k+1)}(\tau, \zeta)-U^{(k)}(\tau, \zeta)+V^{(k)}(\tau, \zeta)|\leq M\tau^2  \left(\dfrac{2}{3}\right)^j.
\end{cases}
\end{align}
\end{l1}
\begin{proof}
We again employ the induction technique to prove this Lemma as well. 

For $n=1$, one has
\begin{align*}
    |U^{(2)}&-U^{(1)}|\leq \displaystyle \bigintssss_{0}^{\tau} \Bigg\{ \dfrac{|U^{(1)}-V^{(1)}|}{2\upsilon}\\
    &+\left(\dfrac{\left(i+\frac{H^2}{\gamma^2}\right)\bigg[ i^2p''(\rho)+\kappa_0^2 n^2(1-p'(\rho))\bigg]\gamma_{w(\upsilon)}^2f(w(\upsilon))}{2i^2(\upsilon)w^2(\upsilon)}+2-\upsilon^2\right)\dfrac{|U^{(1)}-V^{(1)}|\upsilon}{2F}\\
    &+\left(\dfrac{\left(i+\frac{H^2}{\gamma^2}\right)\bigg[ i^2p''(\rho)+\kappa_0^2 n^2(1-p'(\rho))\bigg]\gamma_{w(\upsilon)}^2f(w(\upsilon))}{2i^2(\upsilon)w^2(\upsilon)}+2-\upsilon^2\right)\dfrac{|U^{(1)}|\upsilon}{F}\Bigg\} (\upsilon, \chi_+(\upsilon; \tau, \zeta)) d\upsilon\\
    &\leq \sum_{j=0}^1 \left(\dfrac{2}{3}\right)^j \displaystyle \int_{0}^{\tau} \left(\dfrac{M\upsilon}{2}+2K_\delta M\upsilon^3\right)(\upsilon, \chi_+(\upsilon)) d\upsilon\\
    &\leq M\tau^2\bigg[\dfrac{1}{4}+\dfrac{\delta^2 K_\delta}{2}\bigg]\sum_{j=0}^1 \left(\dfrac{2}{3}\right)^j \leq M\tau^2 \left(\dfrac{2}{3}\right).
\end{align*}
Similarly one can easily prove that $|V^{(2)}-V^{(1)}|\leq M\tau^2\left(\dfrac{2}{3}\right)$ and $|U^{(2)}-V^{(2)}-U^{(1)}+V^{(1)}|\leq M\tau^2\left(\dfrac{2}{3}\right)$. 

We now assume that the Lemma is true for $n=k$, and then we proceed to prove that the Lemma is true for $n=k+1$.

We first use \eqref{eq: 4.24} to obtain
\begin{align*}
   &|V^{(k+1)}(\tau, \zeta)-V^{(k)}(\tau, \zeta)| \\
   &\leq \displaystyle \bigintssss_{0}^{\tau} \bigg\{ \dfrac{|V^{(k)}-U^{(k)}-V^{(k-1)}+U^{(k-1)}|}{2\upsilon}\\
   &+\left(\dfrac{\left(i+\frac{H^2}{\gamma^2}\right)\bigg[ i^2p''(\rho)+\kappa_0^2 n^2(1-p'(\rho))\bigg]\gamma_{w(\upsilon)}^2f(w(\upsilon))}{2i^2(\upsilon)w^2(\upsilon)}+2-\upsilon^2\right)\dfrac{|U^{(k)}-V^{(k)}-U^{(k-1)}+V^{(k-1)}|\upsilon}{2F}\\
    &+\left(\dfrac{\left(i+\frac{H^2}{\gamma^2}\right)\bigg[ i^2p''(\rho)+\kappa_0^2 n^2(1-p'(\rho))\bigg]\gamma_{w}^2f(w)}{2i^2(\upsilon)w^2(\upsilon)}+2(1-\upsilon^2)\right)\dfrac{|V^{(k)}-V^{(k-1)}|\upsilon}{F}\bigg\} (\upsilon, \chi_-(\upsilon; \tau, \zeta)) d\upsilon\\
    & \leq \displaystyle \bigintssss_{0}^{\tau} \bigg\{M\upsilon\left(\dfrac{2}{3}\right)^{k-1}+3K_\delta M\upsilon^3\left(\dfrac{2}{3}\right)^{k-1}\bigg\}d\upsilon=\left(\dfrac{2}{3}\right)^{k-1}\bigg[\dfrac{M\tau^2}{2}+\dfrac{3MK_\delta \tau^4}{4}\bigg]\\
    & \leq M\tau^2 \bigg[\dfrac{1}{2}+\dfrac{3K_\delta \delta^2}{4}\bigg]\left(\dfrac{2}{3}\right)^{k-1}\leq M\tau^2 \left(\dfrac{2}{3}\right)^{k}.
\end{align*}
The estimate of $|U^{(k+1)}-U^{(k)}|$ can be obtained in a similar manner. For the sake of brevity, we omit the details here.

Now we proceed to estimate the quantity $|U^{(k+1)}-V^{(k+1)}-U^{(k)}+V^{(k)}|$. A straightforward calculation leads us to
\begin{align}\label{eq: 4.39}
    |U^{(k+1)}-V^{(k+1)}-U^{(k)}+V^{(k)}|\leq I_3+I_4,
\end{align}
where 
\begin{align*}
    I_3&=\displaystyle \bigintssss_{\tau_1}^{\tau} \bigg\{ \dfrac{|U^{(k)}-V^{(k)}-U^{(k-1)}+V^{(k-1)}|}{2\upsilon}(\upsilon, \chi_+(\upsilon))+\dfrac{|U^{(k)}-V^{(k)}-U^{(k-1)}+V^{(k-1)}|}{2\upsilon}(\upsilon, \chi_-(\upsilon))\\
    &\hspace{-2 cm}+\left(\dfrac{\left(i+\frac{H^2}{\gamma^2}\right)\bigg[ i^2p''(\rho)+\kappa_0^2 n^2(1-p'(\rho))\bigg]\gamma_{w(\upsilon)}^2f(w(\upsilon))}{2i^2(\upsilon)w^2(\upsilon)}+2-\upsilon^2\right)\dfrac{|U^{(k)}-V^{(k)}-U^{(k-1)}+V^{(k-1)}|\upsilon}{2F}(\upsilon, \chi_+(\upsilon))\\
    &\hspace{-2 cm}+\left(\dfrac{\left(i+\frac{H^2}{\gamma^2}\right)\bigg[ i^2p''(\rho)+\kappa_0^2 n^2(1-p'(\rho))\bigg]\gamma_{w(\upsilon)}^2f(w(\upsilon))}{2i^2(\upsilon)w^2(\upsilon)}+2-\upsilon^2\right)\dfrac{|U^{(k)}-V^{(k)}-U^{(k-1)}+V^{(k-1)}|\upsilon}{2F}(\upsilon, \chi_-(\upsilon))\\
    &\hspace{-2 cm}+\left(\dfrac{\left(i+\frac{H^2}{\gamma^2}\right)\bigg[ i^2p''(\rho)+\kappa_0^2 n^2(1-p'(\rho))\bigg]\gamma_{w(\upsilon)}^2f(w(\upsilon))}{2i^2(\upsilon)w^2(\upsilon)}+2(1-\upsilon^2)\right)\dfrac{|U^{(k)}-U^{(k-1)}|\upsilon}{F}(\upsilon, \chi_+(\upsilon)\\
    &\hspace{-2 cm}+\left(\dfrac{\left(i+\frac{H^2}{\gamma^2}\right)\bigg[ i^2p''(\rho)+\kappa_0^2 n^2(1-p'(\rho))\bigg]\gamma_{w(\upsilon)}^2f(w(\upsilon))}{2i^2(\upsilon)w^2(\upsilon)}+2(1-\upsilon^2)\right)\dfrac{|V^{(k)}-V^{(k-1)}|\upsilon}{F}(\upsilon, \chi_-(\upsilon)\bigg\}  d\upsilon
\end{align*}
and 
\begin{align*}
    I_4&=\displaystyle \bigintssss_{0}^{\tau_1} \bigg\{ \dfrac{|V^{(k)}-U^{(k)}-V^{(k-1)}+U^{(k-1)}|}{2\upsilon}\\
    &\hspace{-2 cm}+\left(\dfrac{\left(i+\frac{H^2}{\gamma^2}\right)\bigg[ i^2p''(\rho)+\kappa_0^2 n^2(1-p'(\rho))\bigg]\gamma_{w(\upsilon)}^2f(w(\upsilon))}{2i^2(\upsilon)w^2(\upsilon)}+2-\upsilon^2\right)\dfrac{|U^{(k)}-V^{(k)}-U^{(k-1)}+V^{(k-1)}|\upsilon}{2F}\\
    &\hspace{-2 cm}+\left(\dfrac{\left(i+\frac{H^2}{\gamma^2}\right)\bigg[ i^2p''(\rho)+\kappa_0^2 n^2(1-p'(\rho))\bigg]\gamma_{w(\upsilon)}^2f(w(\upsilon))}{2i^2(\upsilon)w^2(\upsilon)}+2(1-\upsilon^2)\right)\dfrac{|V^{(k)}-V^{(k-1)}|\upsilon}{F}\bigg\} (\upsilon, \chi_+(\upsilon)) d\upsilon.
\end{align*}
Hence by the induction assumption, we have
\begin{align*}
    I_3&\leq \displaystyle \bigintssss_{\tau_1}^{\tau} \bigg\{  M\upsilon\left(\dfrac{2}{3}\right)^{k-1}+6K_\delta M\upsilon^3 \left(\dfrac{2}{3}\right)^{k-1}\bigg\}d\upsilon=\left(\dfrac{2}{3}\right)^{k-1}\bigg\{\dfrac{M}{2}(\tau^2-\tau_1^2)+\dfrac{3}{2}K_\delta M (\tau^4-\tau_1^4)\bigg\},\\
    I_4 &\leq \displaystyle \int_{0}^{\tau} \bigg\{  \dfrac{M\upsilon}{2} \left(\dfrac{2}{3}\right)^{k-1}+3K_\delta M\upsilon^3 \left(\dfrac{2}{3}\right)^{k-1}\bigg\}d\upsilon=\left(\dfrac{2}{3}\right)^{k-1}\bigg\{\dfrac{M}{4}\tau_1^2+\dfrac{3K_\delta M\tau_1^4}{4}\bigg\}.
\end{align*}
which is used in \eqref{eq: 4.39} to obtain
\begin{align*}
    |U^{(k+1)}-V^{(k+1)}-U^{(k)}+V^{(k)}|&\leq \left(\dfrac{2}{3}\right)^{k-1}\bigg\{\dfrac{M}{2}(\tau^2-\tau_1^2)+\dfrac{3}{2}K_\delta M (\tau^4-\tau_1^4)+\dfrac{M}{4}\tau_1^2+\dfrac{3K_\delta M }{4}\tau_1^4\bigg\}\\
    &=\left(\dfrac{2}{3}\right)^{k-1} \bigg\{\dfrac{M\tau^2}{2}+\dfrac{3}{2}K_\delta M \tau^4-\dfrac{M}{4}\tau_1^2-\dfrac{3K_\delta M }{4}\tau_1^4\bigg\}\\
    &\leq M\tau^2 \left(\dfrac{1}{2}+\dfrac{3}{2}\delta^2 K_\delta\right)\left(\dfrac{2}{3}\right)^{k-1}\leq M\tau^2\left(\dfrac{2}{3}\right)^{k}.
\end{align*}
Combining all these estimates proves the lemma.
\end{proof}
\subsection{Existence and uniqueness of the solution in the $(t-r)$ plane}
We now discuss the properties of $(U^{(k)}, V^{(k)})(\tau, \zeta)$ defined by \eqref{eq: 4.24}. It is easy to see that from Lemma \ref{l-4.3} that $(U^{(k)}, V^{(k)})(\tau, \zeta)$ are uniformly convergent sequences of continuous functions which means that the limit function $(U, V)$ must be continuous too and it must satisfy
\begin{align}\label{eq: 4.40}
    |U(\tau, \zeta)|\leq 3M\tau^2, ~~|V(\tau, \zeta)|\leq 3M\tau^2,~~|U(\tau, \zeta)-V(\tau, \zeta)|\leq 3M\tau^2
\end{align}
for all $(\tau, \zeta)\in \Omega$. Further, the functions $(U, V)$ satisfy \eqref{eq: 4.24} and the boundary conditions $(U(0, \zeta))=V(0, \zeta)=0$. Also, in view of \eqref{eq: 4.23} we find that $\zeta=0$ if and only if $\tau=\tau_1$. Applying this in \eqref{eq: 4.24} implies that $U(\tau, 0)=b_1(\tau)$. 

We now verify that the boundary conditions $U_\tau(0, \zeta)=V_\tau(0, \zeta)=0$ are satisfied by $(U, V)$. In view of \eqref{eq: 4.40} and the fundamental theorem of calculus, it is easy to observe that the derivatives of $U(\tau, \zeta)$ and $V(\tau, \zeta)$ with respect to $\tau$ are continuous. We now prove the existence of $(U_\zeta, V_\zeta)$ near $\tau=0$. For the sake of brevity, we only discuss the case $\zeta\geq \bar{\chi}(\tau)$ as the other case can be discussed in a similar manner. Now we differentiate \eqref{eq: 4.24} wrt $\zeta$ to yield the system of integral equations of the form
\begin{align}
\begin{cases}
V_\zeta^{(k)}(\tau, \zeta)&=\displaystyle \bigintssss_0^{\tau} \Bigg\{ \dfrac{V_\chi^{(k-1)}-U_\chi^{(k-1)}}{2\upsilon}+E_{2\chi}\upsilon\\
&\hspace{-1 cm}+\left(\dfrac{\left(i+\frac{H^2}{\gamma^2}\right)\bigg[ i^2p''(\rho)+\kappa_0^2 n^2(1-p'(\rho))\bigg]\gamma_{w(\upsilon)}^2f(w(\upsilon))}{2i^2(\upsilon)w^2(\upsilon)}+2-\upsilon^2\right)\dfrac{(V_\chi^{(k-1)}-U_\chi^{(k-1)})\upsilon }{2F}\\
&\hspace{-1 cm}-\left(\dfrac{\left(i+\frac{H^2}{\gamma^2}\right)\bigg[ i^2p''(\rho)+\kappa_0^2 n^2(1-p'(\rho))\bigg]\gamma_{w}^2f(w)}{2i^2(\upsilon)w^2(\upsilon)}+2(1-\upsilon^2)\right)\dfrac{V_\chi^{(k-1)}\upsilon}{F}\Bigg\}\dfrac{\partial \chi_-}{\partial \zeta}(\upsilon, \chi_-(\upsilon)) d\upsilon,\vspace{0.3 cm}\\
U_\zeta^{(k)}(\tau, \zeta)&=\displaystyle \bigintssss_0^{\tau} \Bigg\{\dfrac{U_\chi^{(k-1)}-V_\chi^{(k-1)}}{2\upsilon}+E_{1\chi}\upsilon\\
&\hspace{-1 cm}+\left(\dfrac{\left(i+\frac{H^2}{\gamma^2}\right)\bigg[ i^2p''(\rho)+\kappa_0^2 n^2(1-p'(\rho))\bigg]\gamma_{w(\upsilon)}^2f(w(\upsilon))}{2i^2(\upsilon)w^2(\upsilon)}+2-\upsilon^2\right)\dfrac{(U_\chi^{(k-1)}-V_\chi^{(k-1)})\upsilon }{2F}\\
&\hspace{-1 cm}-\left(\dfrac{\left(i+\frac{H^2}{\gamma^2}\right)\bigg[ i^2p''(\rho)+\kappa_0^2 n^2(1-p'(\rho))\bigg]\gamma_{w}^2f(w)}{2i^2(\upsilon)w^2(\upsilon)}+2(1-\upsilon^2)\right)\dfrac{U_\chi^{(k-1)}\upsilon}{F(\upsilon)}\Bigg\}\dfrac{\partial \chi_+}{\partial \zeta}(\upsilon, \chi_+(\upsilon)) d\upsilon,
\end{cases} 
\end{align}
where 
\begin{align}
    \dfrac{\partial \chi_{\pm}}{\partial \zeta}= \exp\left(\displaystyle \int_{\tau}^{\upsilon} \dfrac{\partial \chi_{\pm}}{\partial \chi}(s, \chi_\pm(s; \tau, \zeta))ds\right).
\end{align}
Recalling the expressions of $\lambda_\pm$, we find $\dfrac{\partial \chi_{\pm}}{\partial \chi}=0$ and hence  $\dfrac{\partial \chi_{\pm}}{\partial \zeta}=1$. Also by recalling the expressions of $E_1, E_2$ and the regularities of $a_0, a_1$, we also have
\begin{align*}
    &|E_{1\chi}(\upsilon, \chi_+(\upsilon; \tau, \zeta))-E_{2\chi}(\upsilon, \chi_+(\upsilon; \tau, \zeta))|\\
    &\leq \dfrac{2\upsilon \sqrt{1-\upsilon^2}}{F}(|a_0''|+|a_1''|\upsilon)+\dfrac{|a_1'(\chi_+(\upsilon; \tau, \zeta))|+|a_1'(\chi_-(\upsilon; \tau, \zeta))|}{F}\upsilon^3\\
    &+\left(\dfrac{\left(i+\frac{H^2}{\gamma^2}\right)\bigg[ i^2p''(\rho)+\kappa_0^2 n^2(1-p'(\rho))\bigg]\gamma_{w}^2f(w)}{2i^2(\upsilon)w^2(\upsilon)}+2(1-\upsilon^2)\right)\dfrac{|a_0'(\chi_+(\upsilon; \tau, \zeta))-a_0'(\chi_-(\upsilon; \tau, \zeta))|}{F}\\
    &\leq K_\delta M |\chi_+(\upsilon; \tau, \zeta)-\chi_-(\upsilon; \tau, \zeta)|+2K_\delta \upsilon (M+M\upsilon)+2K_\delta M \upsilon^3\\
    & \leq \dfrac{2K_\delta \tau^3}{3}K_\delta M+2K_\delta \tau(M+M\tau)+2K_\delta M \tau^3\\
    & \leq 2K_\delta M\tau\left(\dfrac{\delta^2 K_\delta}{3}+1+\delta+\delta^2\right) \leq K_\delta M\tau.
\end{align*}
One can perform the same iterative process used in Lemma \ref{l-4.2} and Lemma \ref{l-4.3} to conclude the uniform convergence of $(U_\zeta^{(k)}, V_\zeta^{(k)})$ which is enough to prove that the limit functions $(U_\zeta^{(k)}, V_\zeta)$ are continuous and satisfy $U_\zeta(0, \zeta)=V_\zeta(0, \zeta)=0$. Since all partial derivatives of $(U, V)$ are continuous, $(U, V)$ must be differentiable. Also since $(U, V)$ satisfies the boundary conditions \eqref{eq: 4.19}, $(U, V)$ must be a solution of the system \eqref{eq: 4.22} with the boundary data \eqref{eq: 4.19}.

Now we proceed to prove the uniqueness of the solution. On the contrary let us assume that there exist two solutions of the system given by $(U_1, V_1)$ and $(U_2, V_2)$. Then we denote $\Tilde{U}=U_2-U_1$ and $\Tilde{V}=V_2-V_1$ to observe that the functions $(\Tilde{U}, \Tilde{V})$ satisfy the following  integral equations
{
\footnotesize
\begin{align}\label{eq: 4.43}
\begin{cases}
\tilde{V}(\tau, \zeta)&=\displaystyle \bigintssss_0^{\tau} \Bigg\{ \dfrac{\tilde{V}-\Tilde{R}}{2\upsilon}+\Bigg[\dfrac{\left(i+\frac{H^2}{\gamma^2}\right)\bigg[ i^2p''(\rho)+\kappa_0^2 n^2(1-p'(\rho))\bigg]\gamma_{w(\upsilon)}^2f(w(\upsilon))}{2i^2(\upsilon)w^2(\upsilon)}+2-\upsilon^2\Bigg]\dfrac{(\tilde{V}-\tilde{U})\upsilon }{2F}\\
&\hspace{2 cm}-\Bigg[\dfrac{\left(i+\frac{H^2}{\gamma^2}\right)\bigg[ i^2p''(\rho)+\kappa_0^2 n^2(1-p'(\rho))\bigg]\gamma_{w(\upsilon)}^2f(w(\upsilon))}{2i^2(\upsilon)w^2(\upsilon)}+2(1-\upsilon^2)\Bigg]\dfrac{\tilde{V}\upsilon}{F}\Bigg\}(\upsilon, \chi_-(\upsilon)) d\upsilon,\vspace{0.3 cm}\\
\tilde{U}(\tau, \zeta)&=\begin{cases}
    \displaystyle \bigintssss_0^{\tau} \Bigg\{\dfrac{\tilde{U}-\tilde{V}}{2\upsilon}+\Bigg[\dfrac{\left(i+\frac{H^2}{\gamma^2}\right)\bigg[ i^2p''(\rho)+\kappa_0^2 n^2(1-p'(\rho))\bigg]\gamma_{w(\upsilon)}^2f(w(\upsilon))}{2i^2(\upsilon)w^2(\upsilon)}+2-\upsilon^2\Bigg]\dfrac{(\tilde{U}-\tilde{V})\upsilon }{2F}\\
    \hspace{0.5 cm}-\Bigg[\dfrac{\left(i+\frac{H^2}{\gamma^2}\right)\bigg[ i^2p''(\rho)+\kappa_0^2 n^2(1-p'(\rho))\bigg]\gamma_{w(\upsilon)}^2f(w(\upsilon))}{2i^2(\upsilon)w^2(\upsilon)}+2(1-\upsilon^2)\Bigg]\dfrac{\tilde{U}\upsilon}{F(\upsilon)}\Bigg\}(\upsilon, \chi_+(\upsilon)) d\upsilon,\\
    \hspace{12 cm}\zeta\geq \bar{\chi}(\tau),\vspace{0.3 cm}\\
    b_1(\tau_1)+\displaystyle \bigintssss_{\tau_1}^{\tau} \Bigg\{\dfrac{\tilde{U}-\tilde{V}}{2\upsilon}+\Bigg[\dfrac{\left(i+\frac{H^2}{\gamma^2}\right)\bigg[ i^2p''(\rho)+\kappa_0^2 n^2(1-p'(\rho))\bigg]\gamma_{w(\upsilon)}^2f(w(\upsilon))}{2i^2(\upsilon)w^2(\upsilon)}+2-\upsilon^2\Bigg]\dfrac{(\tilde{U}-\tilde{V})\upsilon }{2F}\\
    \hspace{0.5 cm}-\Bigg[\dfrac{\left(i+\frac{H^2}{\gamma^2}\right)\bigg[ i^2p''(\rho)+\kappa_0^2 n^2(1-p'(\rho))\bigg]\gamma_{w(\upsilon)}^2f(w(\upsilon))}{2i^2(\upsilon)w^2(\upsilon)}+2(1-\upsilon^2)\Bigg]\dfrac{\tilde{U}\upsilon}{F(\upsilon)}\Bigg\}(\upsilon, \chi_+(\upsilon)) d\upsilon,\\
    \hspace{12 cm}0\leq \zeta<\bar{\chi}(\tau).
\end{cases} 
\end{cases}
\end{align}
}
It can be observed that $(\tilde{U}, \tilde{V})$ satisfies the inequalities \eqref{eq: 4.38} for $k\geq 1$, therefore a repeated insertion of these inequalities in the right-hand side of \eqref{eq: 4.43} asserts that $|\tilde{U}|\leq \widetilde{M}\left(\dfrac{2}{3}\right)^k$ and $|\tilde{V}|\leq \widetilde{M}\left(\dfrac{2}{3}\right)^k$ for some positive constant $\widetilde{M}>0$ and for any $k$. Therefore, we must have $\tilde{U}=\tilde{V}=0$, which proves the uniqueness of the solution.

Finally, we note by \eqref{eq: 4.13}, \eqref{eq: 4.15} and \eqref{eq: 4.16} that the initial boundary value problem \eqref{eq: 4.5}, \eqref{eq: 4.11} is equivalent to the initial boundary value problem \eqref{eq: 4.19}, \eqref{eq: 4.20}. Therefore, we have the following result in the partial hodograph plane:
\begin{t1}\label{t: 4.1}
Under the assumption \eqref{eq: 4.12}, the boundary value problem \eqref{eq: 4.5} and \eqref{eq: 4.11} have a unique classical solution in the neighbourhood of the point $B'(0, r_2)$ in the region $A'B'C'$.
\end{t1}
\section{Solutions in the physical plane}
We now recover a solution to our main problem in the physical plane by utilizing the solution obtained in Theorem \ref{t: 4.1} with the help of an inverse transformation.
\subsection{Inversion}
In view of \eqref{eq: 4.1} one has
\begin{align}\label{eq: 5.1}
    \dfrac{\partial x}{\partial t}=\dfrac{\theta_y}{J},~~\dfrac{\partial y}{\partial t}=-\dfrac{\theta_x}{J},~~\dfrac{\partial x}{\partial r}=\dfrac{\sin \omega\omega_y}{J},~~\dfrac{\partial y}{\partial r}=-\dfrac{\sin \omega\omega_x}{J},
\end{align}
where $J$ is the Jacobian defined in \eqref{eq: 4.3}. Therefore, a straightforward calculation leads us to 
\begin{align}\label{eq: 5.2}
    \begin{cases}
        \theta_x=t\sin r(W+Z)-\sqrt{1-t^2}\cos r (W-Z), \vspace{0.2 cm}\\
        \theta_y=-t\cos r(W+Z)-\sqrt{1-t^2}\sin r (W-Z), \vspace{0.2 cm}\\
        \varpi_x=-\dfrac{\hat{F}_1(t)}{4i^2(t)w^2(t)t}\left(t\sin r(W-Z) -\sqrt{1-t^2} \cos r (W+Z)\right), \vspace{0.2 cm}\\
        \varpi_y=\dfrac{\hat{F}_1(t)}{4i^2(t)w^2(t)t}\left(t\cos r (W-Z)+\sqrt{1-t^2} \sin r (W+Z)\right).
    \end{cases}
\end{align}
Therefore, \eqref{eq: 5.1} implies
\begin{align}\label{eq: 5.3}
    \begin{cases}
        x_t=-\dfrac{t\cos r (W+Z)+\sqrt{1-t^2} \sin r (W-Z)}{4FWZ}, \vspace{0.2 cm}\\
        y_t=-\dfrac{t\sin r(W+Z)-\sqrt{1-t^2}\cos r (W-Z)}{4FWZ}, \vspace{0.2 cm}\\
        x_r=\dfrac{1}{4t\sqrt{1-t^2}WZ}\left(t\cos r(W-Z) +\sqrt{1-t^2} \sin r (W-Z)\right), \vspace{0.2 cm}\\
        y_r=\dfrac{1}{4t\sqrt{1-t^2}WZ}\left(t\sin r(W-Z) -\sqrt{1-t^2} \cos r (W+Z)\right).
    \end{cases}
\end{align}
Since 
\begin{align*}
    \dfrac{dx(t, r_-(t))}{dt}=x_t+x_r\dfrac{dr_-(t)}{dt}, ~~\dfrac{dy(t, r_-(t))}{dt}=y_t+y_r\dfrac{dr_-(t)}{dt}
\end{align*}
where $\dfrac{dr_-(t)}{dt}=\lambda_-(t)=\dfrac{t^2\sqrt{1-t^2}}{F(t)}$. 

Clearly 
\begin{align}\label{eq: 5.4}
    \begin{cases}
        \dfrac{dx(t, r_-(t))}{dt}=-\dfrac{\sqrt{1-t^2}\sin r+t\cos r}{2F(t)Z(t, r)},\\
        \dfrac{dy(t, r_-(t))}{dt}=-\dfrac{t\sin r-\sqrt{1-t^2}\cos r}{2F(t)Z(t, r)}.
    \end{cases}
\end{align}
Let $\Omega':=\{(t, r)| 0\leq t\leq \delta, \bar{r}(t)-\delta\leq r\leq \bar{r}(t)\}$ be the image of the region $\Omega$ in the $(t-r)$ plane. Then we can use \eqref{eq: 5.4} to define $\hat{x}=x(\hat{t}, \hat{r})$ for any point $(\hat{t}, \hat{r})\in \Omega'$ given by:
\begin{align}
    x(\hat{t}, \hat{r})=\begin{cases}
        \hat{\theta}^{-1}(\tilde{r})-\displaystyle \int_{0}^{\hat{t}}\dfrac{t(\sqrt{1-t^2}\sin r_-(t; \hat{t}, \hat{r}))+t\cos r_-(t; \hat{t}, \hat{r})}{2F(t)Z(t; r_-(t; \hat{t}, \hat{r}))}dt,~~~~~~~\hat{r}\leq r(\hat{t}),\\
        \psi^{-1}(\tilde{\theta}^{-1}(\tilde{r}))-\displaystyle \int_{\tilde{t}}^{\hat{t}}\dfrac{t(\sqrt{1-t^2}\sin r_-(t; \hat{t}, \hat{r}))+t\cos r_-(t; \hat{t}, \hat{r})}{2F(t)Z(t; r_-(t; \hat{t}, \hat{r}))}dt,~~~~~~~\hat{r}> r(\hat{t}),
    \end{cases}
\end{align}
where $\hat{\theta}^{-1}, \tilde{\theta}^{-1}$ and $\psi^{-1}$ denote the inverses of $\hat{\theta}, \tilde{\theta}$ and $\psi$, respectively with
\begin{align}
    r_-(t; \hat{t}, \hat{r})=\hat{r}+\displaystyle \int_{t}^{\hat{t}} \dfrac{t^2\sqrt{1-t^2}}{F(t)}dt, ~~\tilde{r}=\hat{r}+\displaystyle \int_{0}^{\hat{t}} \dfrac{t^2\sqrt{1-t^2}}{F(t)}dt,~~r(\hat{t})=r_2-\displaystyle \int_{0}^{\hat{t}} \dfrac{t^2\sqrt{1-t^2}}{F(t)}dt,
\end{align}
where the numbers $\tilde{t}$ and $\tilde{r}$ can be determined by 
\begin{align}
    \tilde{r}=r_2+\displaystyle \int_{0}^{\tilde{t}} \dfrac{t^2\sqrt{1-t^2}}{F(t)}dt, ~~\tilde{r}=r_2+\displaystyle \int_{\tilde{t}}^{\hat{t}} \dfrac{t^2\sqrt{1-t^2}}{F(t)}dt.
\end{align}
Note that the intersection point of the negative characteristic $r=r_-(t; \hat{t}, \hat{r})$ and positive characteristic $r=r_+(t; 0, r_2)$ is denoted by the point $(\tilde{t}, \tilde{r})$. 

One can define the value $\hat{y}=y(\hat{t}, \hat{r})$ in a similar fashion, which is given by
\begin{align}
    y(\hat{t}, \hat{r})=\begin{cases}
        \varphi(\hat{\theta}^{-1}(\tilde{r}))-\displaystyle \int_{0}^{\hat{t}}\dfrac{t(t\sin r_-(t; \hat{t}, \hat{r})-\sqrt{1-t^2}\cos r_-(t; \hat{t}, \hat{r}))}{2F(t)Z(t; r_-(t; \hat{t}, \hat{r}))}dt,~~~~~~~\hat{r}\leq r(\hat{t}),\\
        \tilde{\theta}^{-1}(\tilde{r})-\displaystyle \int_{\tilde{t}}^{\hat{t}}\dfrac{t(t\sin r_-(t; \hat{t}, \hat{r})-\sqrt{1-t^2}\cos r_-(t; \hat{t}, \hat{r}))}{2F(t)Z(t; r_-(t; \hat{t}, \hat{r}))}dt,~~~~~~~\hat{r}> r(\hat{t}).
    \end{cases}
\end{align}
Noting the expressions of $y_t, y_r, x_t$ and $x_r$, it is clear that 
$$j:=\dfrac{\partial(x, y)}{\partial(t, r)}=\dfrac{t}{4 F(t)W(t, r)Z(t, r))}\neq 0,~\mathrm{for}~ t\in (0, \delta],$$ which shows that the transfromation $(t, r)\longrightarrow (x, y)$ is injective for $t\in (0, \delta]$ and therefore one should have $t=t(x, y)$ and $r=r(x, y)$ so that we can define the functions
\begin{align}\label{eq: 5.78}
    \theta=\theta(t(x,y), r(x, y)),~~\varpi=\sqrt{1-t^2(x, y)},
\end{align}
We claim that the constructed functions $(\theta, \varpi)(x, y)$ is actually a solution to the boundary value problem for the system \eqref{eq: 2.27}.

\subsection{Verification of solutions in the physical plane}
In view of expression of $\bar{\partial}_+$ and \eqref{eq: 5.2}, we have
\begin{align*}
    \bar{\partial}_+\theta&=\cos \alpha\theta_x+\sin \alpha \theta_y\\
    &= \cos \alpha[t\sin r(W+Z)-\sqrt{1-t^2}\cos r (W-Z)]+\sin \alpha[-t\cos r(W+Z)-\sqrt{1-t^2}\sin r (W-Z)]
\end{align*}
and
\begin{align*}
   \dfrac{4i^2(\varpi)w^2(\varpi)\cos \omega}{F_1(\varpi)} \bar{\partial}_+\varpi&=\cos \alpha \varpi_x+\sin \alpha \varpi_y\\
    &=\bigg\{\sin \alpha \left(t\cos r (W-Z)+\sqrt{1-t^2} \sin r (W+Z)\right)\\
    &\hspace{2 cm}-\cos \alpha\left(t\sin r(W-Z) -\sqrt{1-t^2} \cos r (W+Z)\right)\bigg\}.
\end{align*}
Therefore, we have
\begin{align*}
    \bar{\partial}_+\theta+\dfrac{4i^2(\varpi)w^2(\varpi)\cos \omega}{F_1(\omega)}\bar{\partial}_+\varpi&=2Z\left(\sqrt{1-t^2}\cos (\alpha-\theta)-t\sin (\alpha-\theta)\right)=0.
\end{align*}
Clearly, $\theta(x, y)$ and $\varpi(x, y)$ satisfy first equation of \eqref{eq: 2.27}. Similarly, we can show that the functions $\theta(x, y)$ and $\varpi(x, y)$ satisfy the second equation of \eqref{eq: 2.27}, which shows that $\theta(x, y)$ and $\varpi(x, y)$ is the solution of \eqref{eq: 2.27} with boundary data \eqref{eq: 3.2}. 

In view of $\dfrac{dw}{dq}<0$ and $\dfrac{d\gamma}{dq}>0$, we must have $\dfrac{dw}{d\gamma}<0$ which implies $\gamma:=\gamma(w)$. Thus we combine $w=w(\varpi)$ and \eqref{eq: 5.78} together with the relations $u=\dfrac{w\gamma_w \cos \theta}{\gamma \sin \omega}, v=\dfrac{w\gamma_w \sin \theta}{\gamma \sin \omega}$ to define the functions $(w, u, v)(x, y)$ satisfying
\begin{align}
    w(x, y)=w(\varpi(x, y)),~u=\dfrac{w(\varpi(x, y))\gamma_{w(\varpi(x, y))} \cos \theta(x, y)}{\gamma(w(\varpi(x, y))) \varpi(x, y)},~v=\dfrac{w(\varpi(x, y))\gamma_{w(\varpi(x, y))} \sin \theta(x, y)}{\gamma(w(\varpi(x, y))) \varpi(x, y)}
\end{align}
is the solution to our main problem. We summarise the main result of our article in the following theorem:
\begin{t1}
Let $\widehat{AB}: x=\psi(y)$ be a locally concave characteristic curve of relativistic magnetohydrodynamics flow and $\widehat{BC}: y=\varphi(x)$ be a given sonic curve satisfying $\psi''(y_B)<0$ and $\varphi(x_B)=1$. Further, $(\theta, \varpi)|_{\widehat{BC}}=(\hat{\theta}(x), \hat{\varpi}(x))$ and $(\theta, \varpi)|_{\widehat{AB}}=(\tilde{\theta}(x), \tilde{\varpi}(x))$ are given boundary data on these curves such that $\theta'(x_B)<0$ and $(\varphi' \sin \hat{\theta}+\cos \hat{\theta})(x_B)>0$, then the degenerate boundary value problem for the system \eqref{eq: 2.27} admits a unique classical solution $(\theta, \varpi)$ in the angular domain $ABC$ in the neighbourhood of the point $B$.
\end{t1}
\begin{remark}
Note that, in the Newtonian limit $(c\longrightarrow \infty)$ the system \eqref{eq: 1.1} actually reduces to the 2-D steady magnetogasdynamics system which is given by 
\begin{align*}
\begin{cases}
   (\rho u)_x+(\rho v)_y=0,\\
\left(\rho u^2+p+\dfrac{1}{2}\kappa_0^2 \rho^2\right)_x+\left(\rho uv\right)_y=0,\\
\left(\rho uv\right)_x+\left(\rho v^2+p+\dfrac{1}{2}\kappa_0^2 \rho^2\right)_y=0.
\end{cases}
\end{align*}
Therefore, our results extend the results available in the work of Li and Hu \cite{li2019degenerate} to the magnetogasdynamics system as well especially when the pressure law satisfies $p=A\rho^\gamma$. In particular, if we choose $\kappa_0=0$, one can recover the main result from \cite{li2019degenerate}.
\end{remark}
\section*{Acknowledgments}
\textit{The second author (TRS) expresses his gratitude towards SERB, DST, Government of India (Ref. No. CRG/2022/006297) for its financial support through the core research grant.}
\biboptions{sort&compress}
\bibliographystyle{elsarticle-num}
\bibliography{Reference}%
\end{document}